\documentclass[a4paper,12pt]{article}
\usepackage{amsfonts}
\usepackage[xdvi]{graphicx}
\usepackage[dvips]{color}
\usepackage{amsmath,amssymb}

\newcommand{\bm}[1]{\mbox{\boldmath $#1$}}

\newcommand{\C}{{\mathbb C}}

\newcommand{\Z}{{\mathbb Z}}

\makeatletter

\@addtoreset{equation}{section}

\newcounter{def}[section]
\renewcommand{\thedef}{\stepcounter{def}\thesection.\@arabic\c@def }

\begin{document}
\setlength{\baselineskip}{24pt}
\begin{center}
\textbf{\LARGE{Kovalevskaya exponents and the space of initial conditions of a quasi-homogeneous vector field}}
\end{center}

\setlength{\baselineskip}{14pt}

\begin{center}
Institute of Mathematics for Industry, Kyushu University, Fukuoka,
819-0395, Japan

\large{Hayato CHIBA} \footnote{E mail address : chiba@imi.kyushu-u.ac.jp}
\end{center}
\begin{center}
July 6, 2014; Last modified Aug 15 2015
\end{center}

\begin{center}
\textbf{Abstract}
\end{center}

Formal series solutions and the Kovalevskaya exponents of a quasi-homogeneous polynomial system of differential equations are studied 
by means of a weighted projective space and dynamical systems theory.
A necessary and sufficient condition for the series solution to be a convergent Laurent series is given, which improve the well known 
Painlev\'{e} test.
In particular, if a given system has the Painlev\'{e} property,
an algorithm to construct Okamoto's space of initial conditions is given.
The space of initial conditions is obtained by weighted blow-ups of the weighted projective space,
where the weights for the blow-ups are determined by the Kovalevskaya exponents.
The results are applied to the first Painlev\'{e} hierarchy ($2m$-th order first Painlev\'{e} equation).




\section{Introduction}

A system of polynomial differential equations 
\begin{equation}
\frac{dx_i}{dz} = f_i(x_1 , \cdots  ,x_m,z) + g_i(x_1 , \cdots  ,x_m,z), \quad i=1, \cdots  ,m,
\label{1-1}
\end{equation}
is considered, where $(x_1, \cdots ,x_m,z) \in \C^{m+1}$, 
and $f_i$ and $g_i$ satisfy certain conditions on the quasi-homogeneity (see assumptions (A1) to (A3) in Sec.2.1),
for which Kovalevskaya exponents are well defined.
The first, second, fourth Painlev\'{e} equations and the first Painlev\'{e} hierarchy satisfy these conditions.
This system is investigated with the aid of the $m+1$ dimensional weighted projective space $\C P^{m+1}(p_1, \cdots ,p_m, r, s)$ with the positive weight
$(p_1, \cdots ,p_m, r, s)\in \Z^{m+2}_{> 0} $ determined by the quasi-homogeneity of the system.
The space $\C P^{m+1}(p_1, \cdots ,p_m, r, s)$ is decomposed as
\begin{eqnarray*}
\C P^{m+1}(p_1, \cdots ,p_m, r, s) = \C^{m+1}/\Z_s\,\, \cup \,\, \C P^{m}(p_1, \cdots ,p_m, r), \quad (\text{disjoint}).
\end{eqnarray*}
This implies that the space is a compactification of $\C^{m+1}/\Z_s$ obtained by attaching the $m$-dim weighted projective space 
$\C P^{m}(p_1, \cdots ,p_m, r)$ at infinity.
The lift $\C^{m+1} = \{ (x_1, \cdots ,x_m, z)\}$ of the quotient $\C^{m+1}/\Z_s$ is a natural phase space,
on which the system (\ref{1-1}) is given.
The system is also well defined on the quotient $\C^{m+1}/\Z_s$ because it is invariant under the $\Z_s$ action
due to the quasi-homogeneity assumptions.
Then, the system is continuously extended to the codimension one space $\C P^{m}(p_1, \cdots ,p_m, r)$ attached at infinity.
The asymptotic behavior of solutions of the system will be captured by investigating behavior around
the ``infinity set" $\C P^{m}(p_1, \cdots ,p_m, r)$;
The space $\C P^{m+1}(p_1, \cdots ,p_m, r, s)$ gives a 
suitable compactification of the phase space of the system.

A formal series solution of the form
\begin{equation}
x_i(z) = c_i(z-z_0)^{-p_i} +  a_{i,1}(z-z_0)^{-p_i+1}+ a_{i,2}(z-z_0)^{-p_i+2} + \cdots 
\label{1-2}
\end{equation}
will be considered, where $z_0$ is an arbitrary constant (movable singularity), $c_i$ is a constant and $a_{i,n}$ may include $\log (z-z_0)$.
If $a_{i,n}$ is independent of $z$ and the series is convergent, it provides a Laurent series solution.
If coefficients $c_i$ and $a_{i,n}$ include $n$ arbitrary parameters other than $z_0$, it represents an $n+1$-parameter 
family of solutions.
It will be shown that there exists a singularity (in the sense of a foliation defined by integral curves) of the system 
on the ``infinity set" $\C P^{m+1}(p_1, \cdots ,p_m, r)$,
to which the family (\ref{1-2}) of formal series solutions approaches as $z\to z_0$ (Lemma 3.3).
Hence, the asymptotics of (\ref{1-2}) as $z\to z_0$ can be investigated by local analysis around the singularity.
In particular, the normal form theory of dynamical systems will play an important role.
It will be proved in Thm.3.4 that the eigenvalues of the Jacobi matrix at the singularity coincide with the Kovalevskaya exponents,
which implies that the Kovalevskaya exponents are invariant under smooth coordinate transformations.
By combining the weighted projective space $\C P^{m+1}(p_1, \cdots ,p_m, r, s)$, the Kovalevskaya exponents and the normal form theory,
a necessary and sufficient condition for the series  solutions (\ref{1-2}) to be a convergent Laurent series will be given,
which refines the classical Painlev\'{e} test \cite{Abl, Yos2}.
To give the necessary and sufficient condition, it will be shown that the system (\ref{1-1}) has
formal solutions of the form
\begin{equation}
x_i(z) = c_iT^{-p_i}(1 +  \widetilde{h}_i(\alpha _2 T^{\lambda _2}, \cdots ,\alpha _m T^{\lambda_m},z_0T^r, \varepsilon _0 T^s)), \quad T:=z-z_0,
\label{1-3}
\end{equation}
where $\widetilde{h}_i$ is a formal power series in the arguments,
whose coefficients are polynomial in $\log T$, $\alpha _2, \cdots ,\alpha _{m},z_0, \varepsilon _0$ are arbitrary parameters
and $\lambda _2, \cdots , \lambda _m$ are Kovalevskaya exponents other than the trivial exponent $\lambda _1 = -1$ (Lemma 3.6).
Suppose that $\mathrm{Re}(\lambda _i) \leq 0$ for $i=2, \cdots ,k$ and $\mathrm{Re}(\lambda _i) > 0$ for $i=k+1, \cdots ,m$.
The unstable manifold theorem proves that 
\begin{equation}
x_i(z) = c_iT^{-p_i}(1 +  \widetilde{h}_i(0, \cdots ,0,\alpha _{k+1} T^{\lambda_{k+1}}, \cdots , \alpha _m T^{\lambda_m},z_0T^r, \varepsilon _0 T^s)),
\label{1-4}
\end{equation}
is a convergent series.
Further, the normal form theory provides a necessary and sufficient condition for it to be the Laurent series without $\log T$ (Prop.3.5).
One of the necessary condition is that all Kovalevskaya exponents $\lambda _i$ with $\mathrm{Re}(\lambda _i) > 0$
are positive integers, as is well known as the Painlev\'{e} test.

As mentioned, the family of series solutions (\ref{1-2}) tends to the singularity on the ``infinity set" as $z\to z_0$.
If the series is a convergent Laurent series, an algorithm to resolve the singularity by a weighted blow-up will be given.
The weight for the weighted blow-up is determined by the Kovalevskaya exponents.
In particular, if a given system has the Painlev\'{e} property in the sense that any solutions are meromorphic,
our method provides an algorithm to construct the space of initial conditions.
For a polynomial system, a manifold $\mathcal{M}(z)$
is called the  space of initial conditions if any solutions of the system give global holomorphic sections of the 
fiber bundle $\mathcal{P}=\{ (x,z)\, | \, x\in \mathcal{M}(z), z\in \C\}$ over $\C$.
If the system has $n$-types of Laurent series solutions, then the space of initial conditions is obtained by
$n$-times weighted blow-up, which proves that $\mathcal{M}(z)$ is a smooth algebraic variety
obtained by gluing the spaces of the form $\C^m/\Z_{p_j}$ with some integers $p_j$.

In our previous papers \cite{Chi2, Chi3}, weighted projective spaces and dynamical systems theory are applied to
the study of the 2-dim Painlev\'{e} equations (the first to sixth Painlev\'{e} equations).
In particular, it is shown that these equations are linearized by a local analytic transformation
around a movable pole $z = z_0$ (see Prop.3.5), and the spaces of the initial conditions are obtained by the weighted blow-ups.
In the present paper, the previous result is extended to a general quasi-homogeneous system (\ref{1-1}).
In Sec.4, our theory is applied to the first Painlev\'{e} hierarchy, which is a $2m$-dimensional system of equations ($m=1,2,\cdots $).
The $2m$-dimensional first Painlev\'{e} equation has $m$-types of Laurent series solutions.
A complete list of the Kovalevskaya exponents of Laurent series solutions are given (Thm.4.1).
Further, how to construct the space of initial conditions is demonstrated for the 4-dim first Painlev\'{e} equation.


\section{Settings}


\subsection{Kovalevskaya exponent}

Let us consider the system of differential equations
\begin{equation}
\frac{dx_i}{dz} = f_i(x_1 , \cdots  ,x_m, z) + g_i(x_1 , \cdots  ,x_m,z), \quad i=1, \cdots  ,m,
\label{2-1}
\end{equation}
where $f_i$ and $g_i$ are polynomials in $(x_1, \cdots  ,x_m, z) \in \C^{m+1}$.
We suppose that 
\\[0.2cm]
\textbf{(A1)} $(f_1, \cdots  ,f_m)$ is a quasi-homogeneous vector field satisfying
\begin{equation}
f_i(\lambda ^{p_1}x_1 ,\cdots ,\lambda^{p_m}x_m, \lambda ^r z) = \lambda ^{p_i+1}f_i(x_1 , \cdots ,x_m,z)
\label{2-2}
\end{equation}
for any $\lambda \in \C$ and $i=1, \cdots ,m$, where $(p_1, \cdots  ,p_m, r) \in \Z^{m+1}_{>0}$ is a positive weight.
The positive integers $p_i$ and $r$ are called the weighted degrees of $x_i$ and $z$, respectively.

Put $f_i^A (x_1, \cdots ,x_m):= f_i(x_1, \cdots ,x_m, 0)$ and $f_i^N := f_i - f_i^A$
(i.e. $f_i^A$ and $f_i^N$ are autonomous and nonautonomous parts, respectively).
Obviously they satisfy
\begin{eqnarray*}
& & f_i^A(\lambda ^{p_1}x_1 ,\cdots ,\lambda^{p_m}x_m) = \lambda ^{p_i+1}f^A_i(x_1 , \cdots ,x_m), \\
& & f_i^N(\lambda ^{p_1}x_1 ,\cdots ,\lambda^{p_m}x_m, z) = o(\lambda ^{p_i+1}), \quad |\lambda | \to \infty.
\end{eqnarray*}
We assume that $(g_1, \cdots ,g_m)$ is also small with respect to the above weight;
\\[0.2cm]
\textbf{(A2)} Suppose $(g_1, \cdots ,g_m)$ satisfies 
\begin{eqnarray*}
g_i(\lambda ^{p_1}x_1 ,\cdots ,\lambda^{p_m}x_m, \lambda ^r z) = o(\lambda ^{p_i+1}), \quad |\lambda | \to \infty.
\end{eqnarray*}

We also consider the truncated system
\begin{equation}
\frac{dx_i}{dz} = f^A_i(x_1 , \cdots  ,x_m), \quad i=1, \cdots  ,m.
\label{2-3}
\end{equation}
\\
\textbf{Lemma \thedef.} The truncated system is invariant under the scaling
\begin{equation}
(x_1, \cdots ,x_m,z) \mapsto (\lambda ^{p_1}x_1 ,\cdots ,\lambda^{p_m}x_m, \lambda ^{-1}z ).
\label{2-4}
\end{equation}
Further, if the equation
\begin{equation}
-p_ic_i = f_i^A (c_1, \cdots  ,c_m), \quad i=1, \cdots ,m
\label{2-5}
\end{equation}
has a root $(c_1,\cdots ,c_m)\in \C^m$, $x_i(z) = c_i(z-z_0)^{-p_i}$ is an exact solution of the truncated system for any $z_0 \in \C$.
\\[-0.2cm]

The variational equation of $dx_i/dz =  f_i^A(x_1 , \cdots  ,x_m)$ along the solution $x_i(z) = c_i(z-z_0)^{-p_i}$ is given by
\begin{eqnarray*}
\frac{dy_i}{dz} = \sum^m_{k=1}\frac{\partial f_i^A}{\partial x_k}( c_1(z-z_0)^{-p_1},\cdots , c_m(z-z_0)^{-p_m})y_k,
\quad i=1, \cdots ,m.
\end{eqnarray*}
Substituting $y_i = \gamma_i (z-z_0)^{\lambda -p_i}$ with the aid of Eq.(\ref{2-2}) provides
\begin{eqnarray*}
\sum^m_{k=1}\frac{\partial f_i^A}{\partial x_k}(c_1,\cdots ,c_m)\gamma_k + p_i\gamma_i = \lambda \gamma_i.
\end{eqnarray*}
Hence, $\lambda $ is an eigenvalue of the matrix $Df^A(c_1,\cdots ,c_m)+\mathrm{diag} (p_1,\cdots ,p_m)$.
\\[0.2cm]
\textbf{Definition \thedef.} Fix a root $\{ c_i \}_{i=1}^m$ of the equation $-p_ic_i = f_i^A (c_1, \cdots  ,c_m)$.
The matrix 
\begin{equation}
K=\Bigl\{ \frac{\partial f_i^A}{\partial x_j}(c_1,\cdots ,c_m)+ p_i\delta _{ij} \Bigr\}^m_{i,j=1}
\label{2-6}
\end{equation}
and its eigenvalues are called the Kovalevskaya matrix and the Kovalevskaya exponents, respectively, of the system (\ref{2-1})
associated with $\{ c_i\}^m_{i=1}$.
\\[-0.2cm]

Consider a formal series solution of Eq.(\ref{2-1}) of the form
\begin{equation}
x_i = c_i(z-z_0)^{-p_i} +  a_{i,1}(z-z_0)^{-p_i+1}+ a_{i,2}(z-z_0)^{-p_i+2} + \cdots 
\label{2-7}
\end{equation}
Coefficients $a_{i,j}$ are determined by substituting it into Eq.(\ref{2-1}).
The column vector $a_j = (a_{1,j} ,\cdots ,a_{m,j})^T$ satisfies
\begin{eqnarray}
(K-jI)a_j = (\text{a function of $c_i$ and $a_{i,k}$ with $k<j$}).
\label{2-7b}
\end{eqnarray}
If a positive integer $j$ is not an eigenvalue of $K$, $a_j$ is uniquely determined.
If a positive integer $j$ is an eigenvalue of $K$ and (\ref{2-7b}) has no solutions, we have to introduce a logarithmic term $\log (z-z_0)$
into the coefficient $a_j$. In this case, the system (\ref{2-1}) has no Laurent series solution
of the form (\ref{2-7}) with a given $\{ c_i \}_{i=1}^m$.
If a positive integer $j$ is an eigenvalue of $K$ and (\ref{2-7b}) has a solution $a_j$, then $a_j + v$ is also a solution for any 
eigenvectors $v$.
This implies that the series solution (\ref{2-7}) includes a free parameter in $(a_{1,j} ,\cdots ,a_{m,j})$.
If it includes $k-1$ free parameters other than $z_0$, (\ref{2-7}) represents a $k$-parameter family of
Laurent series solutions.
Hence, the classical Painlev\'{e} test \cite{Abl, Yos2} for the necessary condition for the Painlev\'{e} property is stated
as follows;
\\[0.2cm]
\textbf{Classical Painlev\'{e} test.} 
If the system (\ref{2-1}) satisfying (A1) and (A2) has the Painlev\'{e} property in the sense that 
any solutions are meromorphic, then there exist numbers $\{ c_i \}_{i=1}^m$ such that
all Kovalevskaya exponents except for $-1$ (see below) are positive integers,
and the Kovalevskaya matrix is semisimple.
In this case, (\ref{2-7}) represents an $m$-parameter family of Laurent series solutions.
\\

In Prop.3.5, we will give a necessary and sufficient condition for the series (\ref{2-7}) to be a convergent Laurent series.
The next lemmas are well known \cite{Bor, HuYa}.
\\[0.2cm]
\textbf{Lemma \thedef.} (i) $\lambda = -1$ is always a Kovalevskaya exponent with the eigenvector $(-p_1c_1,\cdots ,-p_mc_m)^T$.
\\
(ii) $\lambda =0$ is a Kovalevskaya exponent associated with $\{ c_i\}_{i=1}^m$ if and only if $\{ c_i\}_{i=1}^m$
is not an isolated root of the equation $-p_ic_i = f_i^A (c_1, \cdots  ,c_m)$.
\\[0.2cm]
\textbf{Lemma \thedef.} Consider the Hamiltonian system
\begin{equation}
\frac{dx_i}{dz} = - \frac{\partial H}{\partial y_i},\,\, \frac{dy_i}{dz} = \frac{\partial H}{\partial x_i}, \quad (i=1,\cdots ,m)
\label{2-8}
\end{equation}
with a holomorphic Hamiltonian satisfying
\begin{equation}
H(\lambda ^{p_1}x_1, \lambda ^{q_1}y_1,\cdots ,\lambda ^{p_m}x_m, \lambda ^{q_m}y_m) 
=\lambda ^{h+1} H(x_1,y_1,\cdots ,x_m,y_m).
\label{2-9}
\end{equation}
If $\lambda $ is a Kovalevskaya exponent, so is $\mu$ given by $\lambda +\mu = h$.
\\[-0.2cm]

In what follows, a Kovalevskaya exponent is called a K-exponent for simplicity.
Let us consider the system (\ref{2-1}) with the assumptions (A1) and (A2).
We show that the K-exponents are invariant under a certain class of coordinate transformations.

Consider a holomorphic transformation
\begin{equation}
x_i = \varphi_i(y_1,\cdots ,y_m), \quad (i=1,\cdots ,m),  
\label{2-10}
\end{equation}
which is locally biholomorphic near the point $(x_1, \cdots ,x_m) = (c_1, \cdots ,c_m)$.
The inverse transformation is denoted by $y_i = \psi_i(x_1,\cdots ,x_m)$.
Suppose that $\varphi _i$ satisfies
\begin{equation}
\varphi_i(\lambda ^{q_1}y_1,\cdots ,\lambda ^{q_m}y_m) = \lambda ^{p_i} \varphi_i(y_1,\cdots ,y_m), \quad \lambda \in \C,
\label{2-11}
\end{equation}
with some $(q_1,\cdots ,q_m)\in \Z^m$ for a given $(p_1,\cdots ,p_m)$ in (A1).
It is easy to see that the inverse satisfies
\begin{eqnarray*}
\psi_i(\lambda ^{p_1}x_1,\cdots ,\lambda ^{p_m}x_m) = \lambda ^{q_i} \psi_i(x_1,\cdots ,x_m).
\end{eqnarray*}
By the transformation, Eq.(\ref{2-1}) is brought into the new system
\begin{equation}
\frac{dy_i}{dz} = \sum^m_{j=1}(D\varphi )^{-1}_{ij}f_j(\varphi (y),z)+ \sum^m_{j=1}(D\varphi )^{-1}_{ij}g_j(\varphi (y),z) 
=: F_i(y,z) + G_i(y,z),
\label{2-12}
\end{equation}
where $y = (y_1,\cdots ,y_m)$ and $\varphi = (\varphi _1,\cdots ,\varphi _m)$.
It is straightforward to show that the new system satisfies the conditions (A1) and (A2), in which $(p_1,\cdots ,p_m)$
is replaced by $(q_1,\cdots ,q_m)$.
Hence, the K-exponents of (\ref{2-12}) with the weight $(q_1,\cdots ,q_m)$ are well defined.
\\[0.2cm]
\textbf{Theorem \thedef.} The K-exponents of the system (\ref{2-12}) coincide with those of (\ref{2-1}).
\\[0.2cm]
\textbf{Proof.} Differentiated by $y_k$, Eq.(\ref{2-11}) yields
\begin{equation}
\frac{\partial \varphi _i}{\partial y_k}(\lambda ^{q_1}y_1,\cdots ,\lambda ^{q_m}y_m)
  = \lambda ^{p_i-q_k}\frac{\partial \varphi _i}{\partial y_k}(y_1,\cdots ,y_m).
\label{2-13}
\end{equation}
Differentiating in $\lambda $ and putting $\lambda =1$ for Eqs.(\ref{2-11}) and (\ref{2-13}), we obtain
\begin{eqnarray}
& & \sum^m_{k=1}\frac{\partial \varphi _i}{\partial y_k}q_ky_k = p_i \varphi _i(y_1, \cdots ,y_m), \label{2-14}\\
& & \sum^m_{l=1}\frac{\partial^2 \varphi _i}{\partial y_k \partial y_l}q_ly_l = (p_i-q_k)\frac{\partial \varphi _i}{\partial y_k}(y_1, \cdots ,y_m).
\label{2-15}
\end{eqnarray}
The $(i,k)$-component of the Kovalevskaya matrix $\widetilde{K}$ of (\ref{2-12}) is given by
\begin{eqnarray*}
\widetilde{K}_{ik} = \sum^m_{j=1}\frac{\partial (D\varphi )^{-1}_{ij}}{\partial y_k}f^A_j(\varphi (y))
+ \sum^m_{j=1}(D\varphi )^{-1}_{ij}\sum^m_{l=1} \frac{\partial f^A_j}{\partial y_l} (\varphi (y)) 
  \frac{\partial \varphi _l}{\partial y_k}(y) + q_i \delta _{ik}.
\end{eqnarray*}
By using the equalities (\ref{2-14}),(\ref{2-15}) and $-p_ic_i = f_i^A (c_1, \cdots  ,c_m)$,
we can show that $\widetilde{K}_{ik}$ is rewritten as
\begin{eqnarray*}
\widetilde{K}_{ik} = \sum^m_{j,l=1}(D\varphi )^{-1}_{ij}(p_j \delta _{jl} + (Df^A)_{jl})(D\varphi )_{lk}.
\end{eqnarray*}
This proves that $\widetilde{K}$ is similar to $K$. $\Box$
\\[0.2cm]
\textbf{Proposition \thedef.}
If the system (\ref{2-1}) has a formal series solution (\ref{2-7}),
it is a convergent series on $0 < |z-z_0| < \varepsilon $ for some $\varepsilon >0$.
\\

In Sec.3, a formal series solution (\ref{2-7}) is regarded as an integral curve on an unstable manifold of 
a certain vector field.
Then, Prop.2.6 immediately follows from the unstable manifold theorem, see also Goriely \cite{Gor} for the same result
for autonomous systems.

Next, let us consider the series solution of the form
\begin{equation}
x_i(z) = c_i(z-z_0)^{-q_i} + \sum^\infty_{n=1}a_{i,n}(z-z_0)^{-q_i + n}, \quad (i=1,\cdots ,m).
\label{2-17}
\end{equation}
Note that the order of the leading term is $q_i$, not $p_i$.
\\[0.2cm]
\textbf{Proposition \thedef.} 
If $0\leq q_i < p_i$ and $c_i \neq 0$ for any $i$, then $q_i=0$ for all $i$.
\\[0.2cm]
\textbf{Proposition \thedef.} 
For the system (\ref{2-1}), we further suppose the following condition.
\\
\textbf{(S)} A fixed point of the truncated system is only the origin, i.e,
\begin{equation}
f_i^A(x_1, \cdots ,x_m) = 0 \quad (i=1, \cdots ,m) \Rightarrow (x_1 ,\cdots ,x_m) = (0, \cdots ,0).
\end{equation}
If $q_i > p_i$ for some $i$, then $c_i = 0$.
\\

Prop.2.7 means that if the order of a pole of $x_i(z)$ is smaller than $p_i$ for all $i = 1,\cdots ,m$,
then (\ref{2-17}) should be a local analytic solution.
Prop.2.8 implies that there are no Laurent series solutions $x_i(z)$ whose pole order is larger than $p_i$.
Proofs of Prop.2.7 and 2.8 are given in Appendix B.
Combining three propositions, we have
\\[0.2cm]
\textbf{Theorem \thedef.} 
If the system (\ref{2-1}) satisfies (A1), (A2) and (S), 
any formal series solutions with a singularity at $z=z_0$ are of the form (\ref{2-7}) such that
$(c_1, \cdots ,c_m) \neq (0, \cdots ,0)$, and they are convergent.
\\[0.2cm]
\textbf{Remark \thedef.} 
If the truncated system has a fixed point other than the origin,
then it has a family of fixed points which forms an algebraic curve on $\C^m$ due to the quasi-homogeneity.
If the truncated system is a Hamiltonian system with the Hamiltonian function $H(x_1, \cdots ,x_m)$,
the assumption (S) implies that a singularity of the algebraic variety defined by $\{H=0 \}$ is isolated.
This fact will be essentially used to study a relationship between the Painlev\'{e} equations and singularity theory.
\\

Due to (A1), the truncated system (\ref{2-3}) is invariant under the the $\Z_s$ action
\begin{equation}
(x_1, \cdots ,x_m,z) \mapsto (\omega ^{p_1}x_1, \cdots ,\omega ^{p_m}x_m, \omega ^{r}z), \quad \omega := e^{2\pi i/s},
\label{2-19}
\end{equation}
if $s = r+1$.
For later purpose, we assume that the full system (\ref{2-1}) is also invariant under the same action;
\\[0.2cm]
\textbf{(A3)} The system (\ref{2-1}) is invariant under the $\Z_s$ action (\ref{2-19}) with $s=r+1$.
\\[0.2cm]
\textbf{Example \thedef.}
The first, second and fourth Painlev\'{e} equations in Hamiltonian forms are given by
\begin{eqnarray}
& & (\text{P}_\text{I}) \left\{ \begin{array}{l}
\displaystyle \frac{dx}{dz} = 6y^2 + z \\[0.2cm]
\displaystyle \frac{dy}{dz} = x,  \\
\end{array} \right.
\\
& & (\text{P}_\text{II}) \left\{ \begin{array}{l}
\displaystyle \frac{dx}{dz} = 2y^3 + yz +\alpha  \\[0.2cm]
\displaystyle \frac{dy}{dz} = x,  \\
\end{array} \right.
\\
& & (\text{P}_\text{IV}) \left\{ \begin{array}{l}
\displaystyle \frac{dx}{dz} = -x^2 + 2xy + 2xz + \alpha   \\[0.2cm]
\displaystyle \frac{dy}{dz} = -y^2+2xy-2yz + \beta,  \\
\end{array} \right.
\end{eqnarray}
where $\alpha, \beta \in \C$ are arbitrary parameters.
These systems satisfy the assumptions (A1) to (A3) with the weights
\begin{eqnarray*}
(\text{P}_\text{I}) & & (p_1,p_2,r,s) = (3,2,4,5), \\
(\text{P}_\text{II}) & & (p_1,p_2,r,s) = (2,1,2,3), \\
(\text{P}_\text{IV}) & & (p_1,p_2,r,s) = (1,1,1,2),
\end{eqnarray*}
where $f = (6y^2 + z, x),\, g = (0,0)$ for $(\text{P}_\text{I})$,
$f = (2y^3 + yz, x),\, g= (\alpha ,0)$ for $(\text{P}_\text{II})$ and
$f = ( -x^2 + 2xy + 2xz, -y^2+2xy-2yz),\, g= (\alpha , \beta)$ for $(\text{P}_\text{IV})$.
In Chiba \cite{Chi2}, these systems are investigated by means of 
the weighted projective spaces $\C P^3(p_1,p_2,r,s)$.
One of the purposes in this paper is to extend the previous result to a general system (\ref{2-1}).
For a given weight $(p_1, \cdots ,p_m, r)$ satisfying (A1) to (A3), 
the weighted projective space $\C P^{m+1}(p_1, \cdots , p_m, r, s)$ gives a suitable compactification of $\C^{m+1}$ 
(the space of the dependent variables and the independent variable),
which is effective to investigate the asymptotic behavior of solutions of the system.
\\[0.2cm]
\textbf{Remark \thedef.}
A few remarks are in order.
If $f = (f_1, \cdots ,f_m)$ is independent of $z$ (i.e. $f_i^N = 0$),
we define $r = 0$ and $s=1$.
In this case, we need not assume (A3); the action is trivial.
All results in this paper hold even in this case.
Some of our analysis is still valid even if $(p_1, \cdots ,p_m)$ includes zeros or negative integers.
See \cite{Chi3} for the detail, in which the third, fifth and sixth Painlev\'{e} equations 
are treated.
In general, a weight $(p_1, \cdots ,p_m, r)$ satisfying (A1) to (A3) is not unique.
All possible weights can be calculated through the Newton diagram of the system (\ref{2-1}),
see \cite{Chi2}.
In this paper, we fix one of the weights. 


\subsection{Weighted projective space}

Consider the weighted $\C^*$-action on $\C^{m+2}$ defined by
\begin{eqnarray*}
(x_1, \cdots ,x_{m}, z,\varepsilon ) \mapsto (\lambda^{p_1} x_1, \cdots , \lambda^{p_m} x_m, \lambda^{r} z, \lambda^{s} \varepsilon ), \quad
\lambda \in \C^* := \C\backslash \{ 0\},
\end{eqnarray*}
where the weights $(p_1, \cdots , p_m,r,s)$ are relatively prime positive integers.
The quotient space
\begin{eqnarray*}
\C P^{m+1}(p_1, \cdots , p_m,r,s) := \C^{m+2}\backslash \{ 0\} /\C^*
\end{eqnarray*}
gives an $m+1$ dimensional orbifold called the weighted projective space.

In order to show that a weighted projective space is indeed an orbifold, we will introduce the inhomogeneous coordinates.
For simplicity, we demonstrate it for a three dimensional space $\C P^3(p,q,r,s)$.

The space $\C P^3(p,q,r,s)$ is defined by the equivalence relation on $\C^4\backslash \{ 0\}$
\begin{eqnarray*}
(x,y,z, \varepsilon ) \sim (\lambda^p x, \lambda ^qy, \lambda ^rz, \lambda ^s \varepsilon ).
\end{eqnarray*}
(i) When $x\neq 0$,
\begin{eqnarray*}
(x,y,z,\varepsilon ) \sim (1,\,\,  x^{-q/p}y, \,\, x^{-r/p}z,\,\, x^{-s/p}\varepsilon )
=:(1, Y_1, Z_1, \varepsilon _1).
\end{eqnarray*}
Due to the choice of the branch of $x^{1/p}$, we also obtain
\begin{eqnarray*}
(Y_1, Z_1, \varepsilon _1) \sim
(e^{-2q\pi i /p}Y_1, e^{-2r\pi i /p}Z_1, e^{-2s\pi i /p}\varepsilon _1),
\end{eqnarray*}
by putting $x \mapsto e^{2\pi i }x$.
This implies that the subset of $\C P^3(p,q,r,s)$ such that $x\neq 0$ is homeomorphic to $\C^3 / \Z_p$,
where the $\Z_p$-action is defined as above.

(ii) When $y\neq 0$, 
\begin{eqnarray*}
(x,y,z,\varepsilon ) \sim (y^{-p/q}x,\,\, 1 ,\,\, y^{-r/q}z,\,\, y^{-s/q}\varepsilon )
=:(X_2, 1, Z_2, \varepsilon _2).
\end{eqnarray*}
Because of the choice of the branch of $y^{1/q}$, we obtain
\begin{eqnarray*}
(X_2, Z_2, \varepsilon _2) \sim (e^{-2p\pi i/q}X_2, e^{-2r\pi i/q}Z_2, e^{-2s\pi i/q} \varepsilon _2).
\end{eqnarray*}
Hence, the subset of $\C P^3(p,q,r,s)$ with $y\neq 0$ is homeomorphic to $\C^3 / \Z_q$.

(iii) When $z\neq 0$,
\begin{eqnarray*}
(x,y,z,\varepsilon ) \sim (z^{-p/r}x,\,\, z^{-q/r}y ,\,\, 1,\,\, z^{-s/r}\varepsilon )
=: (X_3, Y_3, 1, \varepsilon _3).
\end{eqnarray*}
Similarly, the subset $\{ z \neq 0\} \subset \C P^3(p,q,r,s)$ is homeomorphic to $\C^3 / \Z_r$.

(iv) When $\varepsilon \neq 0$,
\begin{eqnarray*}
(x,y,z,\varepsilon ) \sim (\varepsilon ^{-p/s}x,\,\, \varepsilon ^{-q/s}y ,\,\, \varepsilon ^{-r/s}z ,\,\, 1)
=:(X_4, Y_4, Z_4, 1).
\end{eqnarray*}
The subset $\{ \varepsilon  \neq 0\} \subset \C P^3(p,q,r,s)$ is homeomorphic to $\C^3 / \Z_s$.

This proves that the orbifold structure of $\C P^3(p,q,r,s)$ is given by
\begin{eqnarray*}
\C P^3(p,q,r,s) = \C^3/\Z_p \,\, \cup \,\, \C^3/\Z_q \,\, \cup \,\, \C^3/\Z_r \,\, \cup \,\, \C^3/\Z_s.
\end{eqnarray*}
The local charts $(Y_1, Z_1, \varepsilon _1)$, $(X_2, Z_2, \varepsilon _2)$,
$(X_3, Y_3, \varepsilon _3)$ and $(X_4, Y_4, Z_4)$ defined above are called inhomogeneous coordinates
as the usual projective space.
Note that they give coordinates on the lift $\C^3$, not on the quotient $\C^3 / \Z_i\,\, (i=p,q,r,s)$.
Therefore, any equations written in these inhomogeneous coordinates should be invariant under
the corresponding $\Z_i$ actions.

The transformations between inhomogeneous coordinates are give by
\begin{eqnarray}
\left\{ \begin{array}{rrrr}
X_4 = & \varepsilon _1^{-p/s} = & X_2\varepsilon _2^{-p/s} =& X_3\varepsilon _3^{-p/s} \\
Y_4 = & Y_1\varepsilon _1^{-q/s} = & \varepsilon _2^{-q/s}=& Y_3\varepsilon _3^{-q/s}\\
Z_4 = & Z_1\varepsilon_1 ^{-r/s} = & Z_2\varepsilon _2^{-r/s}=& \varepsilon _3^{-r/s}.
\end{array} \right.
\label{2-24}
\end{eqnarray}

An extension to the $m+1$ dimensional case $\C P^{m+1}(p_1, \cdots , p_m,r,s)$ is straightforward.
The orbifold structure is characterized by
\begin{eqnarray*}
\C P^{m+1}(p_1, \cdots , p_m,r,s) = \C^{m+1}/\Z_{p_1} \cup \cdots \cup \C^{m+1}/\Z_{p_m}
\cup  \C^{m+1}/\Z_r \cup \C^{m+1}/\Z_s.
\end{eqnarray*}
The inhomogeneous coordinates are defined as above on each chart.
In what follows, we use the notation $(x_1, \cdots , x_m ,z)$ for the inhomogeneous coordinates of the local chart $\C^{m+1}/\Z_s$
because a system of differential equations will be given on this chart.
For example, the transformation between the inhomogeneous coordinates $(x_1, \cdots , x_m ,z)$ on $\C^{m+1}/\Z_s$ and the $j$-th inhomogeneous coordinates
\\
$(X_1, \cdots ,X_{j-1}, X_{j+1},\cdots , X_m, Z, \varepsilon )$ on $\C^{m+1}/\Z_{p_j}$ is give by
\begin{eqnarray}
& & x_i = X_i \varepsilon ^{-p_i/s} \,\, (i\neq j),\,\, x_j=\varepsilon ^{-p_{j}/s},\,\, z = Z \varepsilon ^{-r/s}.
\label{2-25}
\end{eqnarray}
Hence, the subset $\{ \varepsilon =0\}$ on $\C^{m+1}/\Z_{p_j}$ is attached at ``infinity" of the chart $\C^{m+1}/\Z_s$.


\section{A differential equation on a weighted projective space}

Now we give the system of polynomial differential equations
\begin{equation}
\frac{dx_i}{dz} = f_i(x_1 , \cdots  ,x_m,z) + g_i(x_1 , \cdots  ,x_m,z), \quad i=1, \cdots  ,m,
\label{3-1}
\end{equation}
satisfying (A1), (A2) and (A3) on the  $(x_1, \cdots , x_m ,z)$ coordinates of the space $\C P^{m+1}(p_1, \cdots , p_m,r,s)$.
Note that the inhomogeneous coordinates  $(x_1, \cdots , x_m ,z)$ are coordinates for the lift $\C^{m+1}$ of the quotient $\C^{m+1}/\Z_s
\subset \C P^{m+1}(p_1, \cdots , p_m,r,s)$.
Since the system (\ref{3-1}) is invariant under the $\Z_s$ action (\ref{2-19}), it is well defined on the quotient space $\C^{m+1}/\Z_s$.

Let us express it on the $j$-th local chart $\C^{m+1}/\Z_{p_j}$ by the transformation (\ref{2-25}).
Due to the assumptions, we have
\begin{eqnarray*}
& & f_i(X_1\varepsilon ^{-p_1/s},\cdots ,\varepsilon ^{-p_j/s}, \cdots , X_m\varepsilon ^{-p_m/s}, Z\varepsilon ^{-r/s}) \\
&= & \varepsilon ^{-(1+p_i)/s}f_1(X_1,\cdots ,1,\cdots ,X_m,Z), \\
& & g_i(X_1\varepsilon ^{-p_1/s},\cdots ,\varepsilon ^{-p_j/s}, \cdots , X_m\varepsilon ^{-p_m/s}, Z\varepsilon ^{-r/s}) \\
&= & \varepsilon ^{1-(1+p_i)/s} \times \text{(a polynomial of $X_1,\cdots ,X_m,Z, \varepsilon $)}.  
\end{eqnarray*}
With the aid of these equalities, (\ref{3-1}) is written on the $(X_1, \cdots ,X_{j-1}, X_{j+1}, \allowbreak  \cdots  \allowbreak ,X_m,Z,\varepsilon )$
coordinates of $\C^{m+1}/\Z_{p_j}$ as
\begin{equation}
\left\{ \begin{array}{l}
\displaystyle \frac{dX_i}{d\varepsilon } = \frac{1}{s\varepsilon }\left( p_iX_i - p_j \frac{f_i + \varepsilon G_i}{f_j + \varepsilon G_j}\right),
\quad (i=1,\cdots ,m; i\neq j),  \\[0.2cm]
\displaystyle \frac{dZ}{d\varepsilon } = \frac{1}{s\varepsilon }\left( rZ - \frac{p_j \varepsilon }{f_j + \varepsilon G_j}\right),  \\
\end{array} \right.
\label{3-2}
\end{equation}
where $f_i = f_i(X_1,\cdots ,1,\cdots X_m,Z)$ (the unity is substituted at the $j$-th argument)
and $G_i$ is a polynomial in  $(X_1, \cdots ,X_m,Z,\varepsilon )$ determined by $g_i$.
This system is rational and invariant under the $\Z_{p_j}$ action despite the fact that the coordinate transformation (\ref{2-25}) is not rational.
\\[0.2cm]
\textbf{Proposition \thedef.} Give the system (\ref{3-1}) on the $(x_1, \cdots , x_m ,z)$-coordinates
of the local chart $\C^{m+1}/\Z_s$ of $\C P^{m+1}(p_1, \cdots , p_m,r,s)$.
If the system satisfies the assumptions (A1) to (A3), it induces a well defined rational differential equations 
on $\C P^{m+1}(p_1, \cdots , p_m,r,s)$.
\\[0.2cm]
\textbf{Example \thedef.}
We give the first Painlev\'{e} equation $x' = 6y^2 +z,\, y' = x$ on the fourth local chart $(x,y,z)$ of $\C P^3(3,2,4,5)$.
By (\ref{2-24}), it is transformed into the following equations 
\begin{eqnarray*}
& & \frac{dY_1}{d\varepsilon _1} = \frac{3 - 12Y_1^3 - 2Y_1Z_1}{\varepsilon _1 (-30 Y_1^2 - 5Z_1)},
\quad \frac{dZ_1}{d\varepsilon _1} = \frac{3\varepsilon _1 - 24Y_1^2 Z_1 - 4 Z_1^2}{\varepsilon _1 (-30 Y_1^2 - 5Z_1)}, \label{P11} \\[0.2cm]
& & \frac{dX_2}{d\varepsilon _2} = \frac{-12 - 2Z_2 + 3X_2^2}{5X_2\varepsilon _2}, \quad
\frac{dZ_2}{d\varepsilon _2} = \frac{-2\varepsilon _2  + 4 X_2Z_2}{5X_2\varepsilon _2}, \label{P12} \\[0.2cm]
& & \frac{dX_3}{d\varepsilon _3} = \frac{24Y_3^2 + 4 - 3X_3\varepsilon _3}{-5\varepsilon _3^2}, \quad
\frac{dY_3}{d\varepsilon _3} = \frac{4X_3 - 2Y_3\varepsilon _3}{-5\varepsilon _3^2}, \label{P13}
\end{eqnarray*}
on the other inhomogeneous coordinates.
Although the transformations (\ref{2-24}) have branches, the above equations are rational
because the first  Painlev\'{e} equation satisfies (A1) to (A3) with $(p, q, r, s) = (3,2,4,5)$.
Hence, they define a rational ODE on $\C P^3(3,2,4,5)$ in the sense of an orbifold.
\\[-0.2cm]

It is convenient to rewrite (\ref{3-2}) as an autonomous vector field of the form
\begin{equation}
\left\{ \begin{array}{l}
\displaystyle \frac{dX_i}{dt} = p_iX_i - p_j \frac{f_i + \varepsilon G_i}{f_j + \varepsilon G_j}
\quad (i=1,\cdots ,m; i\neq j),  \\[0.2cm]
\displaystyle \frac{dZ}{dt} =rZ - \frac{p_j \varepsilon }{f_j + \varepsilon G_j}, \\[0.2cm]
\displaystyle \frac{d\varepsilon }{dt} = s\varepsilon.
\end{array} \right.
\label{3-3}
\end{equation}
The new independent variable $t$ parameterizes integral curves of (\ref{3-2}).
Note that how to rewrite the system (\ref{3-2}) as an autonomous vector field is not unique.
To construct the space of initial conditions, rewriting as a polynomial vector field 
\begin{equation}
\left\{ \begin{array}{l}
\displaystyle \frac{dX_i}{dt} = p_iX_i (f_j + \varepsilon G_j) - p_j (f_i + \varepsilon G_i ),
\quad (i=1,\cdots ,m; i\neq j),  \\[0.2cm]
\displaystyle \frac{dZ}{dt} =rZ (f_j + \varepsilon G_j) -p_j \varepsilon, \\[0.2cm]
\displaystyle \frac{d\varepsilon }{dt} = s\varepsilon (f_j + \varepsilon G_j),
\end{array} \right.
\label{3-3c}
\end{equation}
may be more convenient, though in this section we will use the form (\ref{3-3}).
\\[-0.2cm]

Let us investigate the K-exponents of the system (\ref{3-1}).
Let $(c_1 ,\cdots ,c_m)$ be one of the roots of the equation $-p_ic_i = f_i^A (c_1 ,\cdots ,c_m)$,
and consider a series solution (\ref{2-7}).
If it is not a local holomorphic solution, $(c_1, \cdots ,c_m) \neq (0, \cdots ,0)$ due to Prop.2.7.
Assume $c_j \neq 0$.
On the $(X_1, \cdots  ,X_m, Z, \varepsilon )$ coordinates of $\C^{m+1}/\Z_{p_j}$, we obtain
\begin{eqnarray*}
\left\{ \begin{array}{l}
X_i = x_ix_j^{-p_i/p_j} = c_ic_j^{-p_i/p_j}(1 + O(z-z_0)),  \\
Z = zx_j^{-r/p_j} = z c_j^{-r/p_j}(z-z_0)^r(1 + O(z-z_0)), \\
\varepsilon  = x_j^{-s/p_j} = c_j^{-s/p_j}(z-z_0)^s(1 + O(z-z_0)). 
\end{array} \right.
\end{eqnarray*}
In particular, 
\begin{eqnarray*}
X_i \to  c_ic_j^{-p_i/p_j}\,\, (i\neq j), \quad Z,\varepsilon \to 0,
\end{eqnarray*}
as $z\to z_0$.
\\[0.2cm]
\textbf{Lemma \thedef.} The point 
\begin{equation}
(X_1, \cdots  ,X_m, Z, \varepsilon ) = ( c_1c_j^{-p_1/p_j}, \cdots  , c_mc_j^{-p_m/p_j}, 0, 0)
\label{3-5}
\end{equation}
is a fixed point of the vector field (\ref{3-3}).
\\[0.2cm]
\textbf{Proof.}
A fixed point of (\ref{3-3}) satisfying $Z = \varepsilon =0$ is given as a root of 
\begin{eqnarray*}
p_iX_i f_j^A(X_1, \cdots ,1,\cdots ,X_m) - p_jf_i^A(X_1, \cdots ,1,\cdots ,X_m) = 0,\,\, (i\neq j).
\end{eqnarray*}
At the point (\ref{3-5}), the left hand side is estimated as
\begin{eqnarray*}
& & p_ic_ic_j^{-p_i/p_j}f_j^A( c_1c_j^{-p_1/p_j}, \cdots  , c_mc_j^{-p_m/p_j})-p_jf_i^A( c_1c_j^{-p_1/p_j}, \cdots  , c_mc_j^{-p_m/p_j}) \\
&=& p_ic_ic_j^{-p_i/p_j}c_j^{-(1+p_j)/p_j} f_j^A( c_1, \cdots  , c_m) - p_jc_j^{-(1+p_i)/p_j} f_i^A( c_1, \cdots  , c_m).
\end{eqnarray*}
Then, the equality $-p_ic_i = f_i^A (c_1 ,\cdots ,c_m)$ proves that the above quantity actually becomes zero. $\Box$
\\[-0.2cm]

This lemma suggests that the behavior of the series solution (\ref{2-7}) as $z\to z_0$ is governed by 
local properties of the fixed point (\ref{3-5}).
In particular, the dynamical systems theory will be applied to the study of the fixed point.
Due to the orbifold structure, the inhomogeneous coordinates $(X_1, \cdots ,X_m,Z,\varepsilon )$ should be divided by the $\Z_{p_j}$-action (Sec.2.3).
Hence, all points expressed as (\ref{3-5}) obtained by different choices of roots $c_j^{-1/p_j}$ 
represent the same point on the quotient space $\C^{m+1}/\Z_{p_j}$. 


\subsection{Kovalevskaya exponents}

The main theorem in this section is stated as follows.
\\[0.2cm]
\textbf{Theorem \thedef.} The eigenvalues of the Jacobi matrix of the vector field (\ref{3-3}) at the fixed point (\ref{3-5})
are given by $r, s$ and K-exponents of the system (\ref{3-1}) other than the trivial exponent $-1$.
(If we use the polynomial vector field (\ref{3-3c}) instead of (\ref{3-3}), eigenvalues change by a constant factor).
\\[0.2cm]
\textbf{Proof.} We assume $j=m$ for simplicity.
Thus, (\ref{3-3}) is an equation for $(X_1,\cdots , \allowbreak  X_{m-1},Z,\varepsilon )$.

Let $v = (-p_1c_1 , \cdots  ,-p_mc_m)^T$ be the eigenvector of the K-matrix associated with the eigenvalue $-1$ (see Lemma 2.3).
Put $v_1 = (-p_1c_1 , \cdots  ,-p_{m-1}c_{m-1})^T, v_2 = -p_mc_m$ and
\begin{eqnarray*}
P = \left(
\begin{array}{@{\,}cc@{\,}}
I & v_1 \\
0 & v_2
\end{array}
\right),
\end{eqnarray*}
where $I$ denotes the $(m-1)\times (m-1)$ identity matrix.
We obtain
\begin{eqnarray*}
P^{-1}KP = \left(
\begin{array}{@{\,}cc@{\,}}
\widetilde{K} & 0 \\
* & -1
\end{array}
\right),
\end{eqnarray*}
where $\widetilde{K}$ is an $(m-1)\times (m-1)$ matrix whose eigenvalues are K-exponents other than $-1$.
The $(i,j)$ component of $\widetilde{K}$ is given by
\begin{eqnarray*}
\widetilde{K}_{i,j}
= p_i\delta _{ij}+ \frac{\partial f_i^A}{\partial x_j}(c_1, \cdots ,c_m)-\frac{p_ic_i}{p_mc_m}\frac{\partial f_m^A}{\partial x_j}(c_1, \cdots ,c_m).
\end{eqnarray*}

On the other hand, the Jacobi matrix of (\ref{3-3}) at the fixed point is of the form
\begin{eqnarray}
J = \left(
\begin{array}{@{\,}c|cc@{\,}}
\widetilde{J} & * & * \\ \hline
0 & r & * \\
0 & 0 & s  
\end{array}
\right),
\label{J}
\end{eqnarray}
where $\widetilde{J}$ is an $(m-1)\times (m-1)$ matrix whose $(i,j)$ component is given by
\begin{equation}
\widetilde{J}_{i,j} = 
p_i\delta _{ij} - \frac{p_m}{f_m^A}\frac{\partial f_i^A}{\partial x_j} +p_m\frac{f_i^A}{(f_m^A)^2}\frac{\partial f_m^A}{\partial x_j},
\label{3-7b}
\end{equation}
where $f_i^A$ is estimated at the point $( c_1c_m^{-p_1/p_m}, \cdots  , c_{m-1}c_m^{-p_{m-1}/p_m},1)$.
Due to the quasi-homogeneity and the equality $-p_ic_i = f_i^A(c_1, \cdots  ,c_m)$, we have
\begin{eqnarray*}
& & f_i^A( c_1c_m^{-p_1/p_m}, \cdots  , c_{m-1}c_m^{-p_{m-1}/p_m},1)=c_m^{-(1+p_i)/p_m} f_i^A(c_1,\cdots ,c_m) = -p_ic_ic_m^{-(1+p_i)/p_m}, \\
& & \frac{\partial f_i^A}{\partial x_j}( c_1c_m^{-p_1/p_m}, \cdots  , c_{m-1}c_m^{-p_{m-1}/p_m},1)
 = c_m^{-(p_i+1-p_j)/p_m}  \frac{\partial f_i^A}{\partial x_j}( c_1,\cdots  , c_m).
\end{eqnarray*}
By a suitable scaling of the independent variable $z$ of the original system (\ref{3-1}), we can assume without loss of generality that $c_m=1$.
Substituting the above equalities into (\ref{3-7b}) with $c_m=1$, we obtain $\widetilde{J}_{i,j} = \widetilde{K}_{i,j}$. $\Box$
\\[-0.2cm]

Since eigenvalues of the Jacobi matrix are invariant under the actions of diffeomorphisms, 
the K-exponents are invariant under a wide class of coordinate transformations.
In particular, the set of all K-exponents associated with all roots $\{ c_i \}^m_{i=1}$ are invariant
under the automorphisms on $\C P^{m+1}(p_1, \cdots ,p_m,r,s)$.


\subsection{Extended Painlev\'{e} test}

Let $\lambda _2, \cdots , \lambda _m, r, s$ be eigenvalues of the Jacobi matrix $J$ of (\ref{3-3}) at the fixed point (\ref{3-5}),
among which $\lambda _2 ,\cdots , \lambda _m$ are K-exponents of (\ref{3-1}) ($\lambda _1 = -1$ is used for the trivial one).
Put $\hat{X}_i = X_i - c_ic_j^{-p_i/p_j}$.
Then, Eq.(\ref{3-3}) is rewritten as
\begin{equation}
\frac{d\bm{\hat{X}}}{dt} = J\hat{\bm{X}} + F(\bm{\hat{X}}),
 \quad \hat{\bm{X}} = (\hat{X}_1,\cdots ,\hat{X}_{j-1}, \hat{X}_{j+1}, \cdots ,\hat{X}_m, Z, \varepsilon ),
\label{3-7}
\end{equation}
with the nonlinearity $F$.
Suppose that $\mathrm{Re}(\lambda _i)\leq 0$ for $i=2,\cdots ,k$, and $\mathrm{Re}(\lambda _i)> 0$ for $i=k+1,\cdots ,m$.
Since $r$ and $s$ are positive integers, $J$ has exactly $m-k+2$ eigenvalues with positive real parts.
Thus, the system (\ref{3-7}) has an $m-k+2$ dimensional unstable manifold at the origin;
by a suitable linear transformation, (\ref{3-7}) is rewritten as 
\begin{eqnarray*}
\left\{ \begin{array}{ll}
\displaystyle \frac{d\bm{X}_u}{dt} = J_u \bm{X}_u + F_u(\bm{X}_u, \bm{X}_s), & \bm{X}_u\in \C^{m-k+2},  \\[0.2cm]
\displaystyle \frac{d\bm{X}_s}{dt} = J_s \bm{X}_s + F_s(\bm{X}_u, \bm{X}_s), & \bm{X}_s\in \C^{k-1},  \\
\end{array} \right.
\end{eqnarray*}
where real parts of eigenvalues of $J_u$ and $J_s$ are positive and nonpositive, respectively.
Due to the unstable manifold theorem, there exists a local analytic function $\varphi $ satisfying $\varphi (0)=D\varphi (0) = 0$
such that the set $(\bm{X}_u, \varphi (\bm{X}_u))$ expresses the unstable manifold.
Then, 
\begin{equation}
\frac{d\bm{X}_u}{dt} = J_u \bm{X}_u + F_u(\bm{X}_u,  \varphi (\bm{X}_u))
\label{3-8}
\end{equation}
gives the dynamics on the unstable manifold.
The purpose in this section is to prove
\\[0.2cm]
\textbf{Proposition \thedef.}
The system (\ref{3-1}) has an $m-k+1$ parameter family of convergent Laurent series solutions of the form
\begin{eqnarray*}
x_i(z) = c_i(z-z_0)^{-p_i} + \sum^\infty_{n=1}a_{i,n} (z-z_0)^{-p_i+n}, \quad i=1,\cdots ,m, 
\end{eqnarray*}
where $\{ a_{i,n}\}$ includes $m-k$ free parameters other than $z_0$, if and only if
\\
\textbf{(i)} $\lambda _{k+1},\cdots ,\lambda _m$ are positive integers,
\\
\textbf{(ii)} $J_u$ is semi-simple, and
\\
\textbf{(iii)} the system (\ref{3-8}) on the unstable manifold is linearizable by a local analytic transformation.

In particular, the system (\ref{3-1}) has an $m$ parameter family of Laurent series solutions if and only if
\\
\textbf{(i)} all K-exponents except for $\lambda _1 = -1$ are positive integers (classical Painlev\'{e} test),
\\
\textbf{(ii)} the Jacobi matrix $J$ is semisimple, and
\\
\textbf{(iii)} the system (\ref{3-7}) is linearizable by a local analytic transformation.
\\[-0.2cm]

A similar result is also obtained by Goriely \cite{Gor} for autonomous systems.
A linearization of the system (\ref{3-8}) is achieved by Poincar\'{e}-Dulac normal form theory by finite steps.
For the convenience of the reader, a brief review of the normal form theory is given in Appendix A.
In \cite{Chi2,Chi3}, it is proved that the first to sixth Painlev\'{e} equations are linearizable.
\\[0.2cm]
\textbf{Proof of Prop.3.5.}
By the standard perturbation method (the variation of constants method), a general solution of (\ref{3-7}) is expressed as
\begin{eqnarray*}
& & X_i(t) = c_ic_j^{-p_i/p_j} + h_i(\alpha _2 e^{\lambda _2 t},\cdots ,\alpha _m e^{\lambda _m t}, z_0 e^{rt}, \varepsilon _0 e^{st}),
\quad (i=1,\cdots ,m;\, i\neq j), \\
& & Z(t) = z_0e^{rt} + c_j^{1/p_j}\varepsilon _0e^{st} + 
           e^{st} h_{m+1}(\alpha _2 e^{\lambda _2 t},\cdots ,\alpha _m e^{\lambda _m t}, z_0 e^{rt}, \varepsilon _0 e^{st}),
\end{eqnarray*}
and $\varepsilon (t) = \varepsilon _0 e^{st}$, where $\alpha _2,\cdots ,\alpha _m, z_0, \varepsilon _0$ are free parameters determined by
an initial condition.
The functions $h_1, \cdots ,h_{m+1}$ are \textit{formal} power series with $h_i(0) = 0$ whose coefficients are polynomial in $t$. 
More precisely, they are expressed as
\begin{eqnarray*}
h_i(\alpha _2 e^{\lambda _2},\cdots ,\varepsilon _0 e^{st})
= \sum^\infty_{|n|=1} h_{i,n}(t)e^{\langle \lambda, n  \rangle t}, \quad (i=1,\cdots ,m+1; \, i\neq j),
\end{eqnarray*}
where $n = (n_2, \cdots ,n_{m+2}),\, |n| = n_2 + \cdots +n_{m+2}$ and 
$\langle \lambda, n  \rangle  = \lambda _2 n_2 + \cdots + \lambda _m n_m + rn_{m+1} + sn_{m+2}$.
The function $h_{i,n} (t)$ is polynomial in $t$.
Moving to the original coordinates $(x_1, \cdots ,x_m, z)$, we obtain the next lemma,
which can be also proved directly from Eq.(\ref{3-1}).
\\[0.2cm]
\textbf{Lemma \thedef.} The system (\ref{3-1}) has a \textit{formal} series solution of the form
\begin{equation}
x_i(z) = c_iT^{-p_i}(1 +  \widetilde{h}_i(\alpha _2 T^{\lambda _2}, \cdots ,\alpha _m T^{\lambda_m},z_0T^r, \varepsilon _0 T^s)), \quad T:=z-z_0,
\label{3-10}
\end{equation}
where $\widetilde{h}_i$ is a formal power series in the arguments,
whose coefficients are polynomial in $\log T$, and $\alpha _2, \cdots ,\alpha _{m},z_0, \varepsilon _0$ are free parameters.
\\[-0.2cm]

A solution on the unstable manifold is obtained by putting $\alpha _2 = \cdots = \alpha _k=0$;
\begin{eqnarray*}
& & X_i(t) = c_ic_j^{-p_i/p_j} + h_i(0,\cdots ,0, \alpha _{k+1} e^{\lambda _{k+1} t},\cdots ,\alpha _m e^{\lambda _m t}, z_0 e^{rt}, \varepsilon _0 e^{st}), \\
& & Z(t) = z_0e^{rt} + c_j^{1/p_j}\varepsilon _0e^{st} + 
           e^{st} h_{m+1}(0,\cdots ,0, \alpha _{k+1} e^{\lambda _{k+1} t},\cdots ,\alpha _m e^{\lambda _m t}, z_0 e^{rt}, \varepsilon _0 e^{st}),
\end{eqnarray*}
This is a \textit{convergent} series solution because of the unstable manifold theorem.
Moving to the original coordinates $(x_1, \cdots ,x_m, z)$, we obtain
\begin{equation}
x_i(z) = c_iT^{-p_i}(1 +  \widetilde{h}_i(0,\cdots ,0, \alpha _{k+1} T^{\lambda _{k+1}}, \cdots ,\alpha _m T^{\lambda_m},z_0T^r, \varepsilon _0 T^s)), 
\label{3-11}
\end{equation}
where the right hand side is a convergent power series in the arguments,
whose coefficients are polynomial in $\log T$, and $\alpha _{k+1}, \cdots ,\alpha _{m},z_0, \varepsilon _0$ are $m-k+2$ free parameters.
On the parameter space, there are curves $(\alpha _{k+1}(t), \cdots , \alpha _{m}(t), z_0(t), \varepsilon _0(t))$,
on which the above solution represents the same solution.
Hence, (\ref{3-11}) defines an $m-k+1$ parameter family of solutions.
It does not include $\log T$ if and only if the coefficients of $h_i$ do not include polynomial in $t$.
Then, Prop.3.5 immediately follows from Prop.A.3. $\Box$


\subsection{The space of initial conditions}

In this section, we give an algorithm to construct the space of initial conditions for a differential equation having the Painlev\'{e} property.
For a polynomial system, a manifold $\mathcal{M}(z)$
is called the  space of initial conditions if any solutions of the system give global holomorphic sections of the 
fiber bundle $\mathcal{P}=\{ (x,z)\, | \, x\in \mathcal{M}(z), z\in \C\}$ over $\C$.
Okamoto \cite{Oka} constructed the spaces of initial conditions for the first to the sixth Painlev\'{e} equations 
by blow-ups of a Hirzebruch surface eight times and by removing a certain divisor called vertical leaves.
In Chiba \cite{Chi2}, the spaces of initial conditions for the first, second and fourth Painlev\'{e} equations, respectively,
are obtained only by one, two and three times blow-ups with the aid of the weighted projective spaces.
The purpose of this section is to extend this result.
If a given equation having the Painlev\'{e} property has $n$-types Laurent series solutions (i.e. the equation $-p_ic_i = f_i^A(c_1, \cdots ,c_m)$
to determine the leading coefficients has $n$ roots, and the corresponding series solutions are convergent Laurent series),
then the space of initial conditions is obtained by $n$ times weighted blow-up of the weighted projective space.

Suppose that the system (\ref{3-1}) satisfying (A1) to (A3) is given.
Determine the leading coefficients $(c_1, \cdots ,c_m)$ of the formal series solution (\ref{2-7}) 
by solving $-p_ic_i = f_i^A(c_1, \cdots ,c_m)$.
For each $(c_1, \cdots ,c_m)$, we suppose that (\ref{2-7}) is a convergent Laurent series (without $\log (z-z_0)$),
so that the conditions (i) to (iii) of Prop.3.5 are satisfied.
In particular, $\lambda _2, \cdots ,\lambda _k$ satisfy $\mathrm{Re}(\lambda _i) \geq 0$ and $\lambda _{k+1}, \cdots ,\lambda _m$ are positive integers.
Then, the fixed point (\ref{3-5}) is a singularity of the foliation defined by integral curves;
any Laurent series solutions pass through this point.
For each fixed point, we will perform the resolution of singularities.
The procedure is divided into five steps as follows;
\\

\textbf{Step 1.} Due to Prop.2.7, $(c_1, \cdots ,c_m)\neq (0, \cdots ,0)$.
Assume $c_j \neq 0$.
Move to the inhomogeneous coordinates on $\C^{m+1}/\Z_{p_j}$ by (\ref{2-25}) to obtain (\ref{3-3});
\begin{equation}
\left\{ \begin{array}{l}
\displaystyle \frac{dX_i}{dt} = p_iX_i - p_j \frac{f_i + \varepsilon G_i}{f_j + \varepsilon G_j}
\quad (i=1,\cdots ,m; i\neq j),  \\[0.2cm]
\displaystyle \frac{dZ}{dt} =rZ - \frac{p_j \varepsilon }{f_j + \varepsilon G_j}, \\[0.2cm]
\displaystyle \frac{d\varepsilon }{dt} = s\varepsilon.
\end{array} \right.
\label{3-3b}
\end{equation}
Our procedure below is independent of how to rewrite the system (\ref{3-2}) to an autonomous vector field.
For example, it may be more convenient to use the polynomial system (\ref{3-3c}) instead of (\ref{3-3}) when
calculating the normal form at Step 3.

This vector field has a fixed point (singularity) (\ref{3-5}).
In order for the point to be the origin, put $\hat{X}_i = X_i - c_ic_j^{-p_i/p_j}$, which results in
\begin{equation}
\left\{ \begin{array}{l}
\displaystyle \frac{d\bm{\hat{X}}}{dt} = \widetilde{J}\bm{\hat{X}}+\bm{F}(\bm{\hat{X}},Z,\varepsilon ), 
   \quad \bm{\hat{X}} = (\hat{X}_1, \cdots ,\hat{X}_{j-1},\hat{X}_{j+1},\cdots , \hat{X}_m),  \\[0.2cm]
\displaystyle \frac{dZ}{dt} = rZ + \varepsilon F_{m+1}(\bm{\hat{X}},Z,\varepsilon ),  \\[0.2cm]
\displaystyle \frac{d\varepsilon }{dt} = s\varepsilon,
\end{array} \right.
\label{3-12}
\end{equation}
where $\bm{F}$ is a nonlinearity and $ F_{m+1}=-p_j/(f_j + \varepsilon G_j)$.
The matrix $\widetilde{J}$ is a submatrix of $J$, whose eigenvalues are nontrivial K-exponents
$\lambda _2, \cdots , \lambda _m$  (see Eq.(\ref{J})).

\textbf{Step 2.} If the Jacobi matrix $\widetilde{J}$ has eigenvalues 
$\lambda _2, \cdots ,\lambda _k$ having nonpositive real parts, transform (\ref{3-12}) as
\begin{equation}
\left\{ \begin{array}{l}
\displaystyle \frac{d\bm{X}_s}{dt} = J_s\bm{X}_s + \bm{F}_s (\bm{X}_s, \bm{X}_u, Z,\varepsilon ), \quad \bm{X}_s\in \C^{k-1},  \\[0.2cm]
\displaystyle \frac{d\bm{X}_u}{dt} = J_u\bm{X}_u + \bm{F}_u (\bm{X}_s, \bm{X}_u, Z,\varepsilon ),\quad \bm{X}_u\in \C^{m-k},   \\[0.2cm]
\displaystyle \frac{dZ}{dt} = rZ + \varepsilon F_{m+1}(\bm{X}_s, \bm{X}_u,Z,\varepsilon ),  \\[0.2cm]
\displaystyle \frac{d\varepsilon }{dt} = s\varepsilon,
\end{array} \right.
\label{3-13}
\end{equation}
by a linear transformation of $(\hat{X}_1, \cdots , \hat{X}_m)$ (we need not change $Z$ and $\varepsilon $),
where real parts of eigenvalues of $J_u$ and $J_s$ are positive and nonpositive, respectively.
We suppose that $J_u$ is of the Jordan normal form.
Because of Prop.3.5, it is semi-simple; $J_u = \mathrm{diag} (\lambda _{k+1}, \cdots ,\lambda _m ) $.
Due to the unstable manifold theorem, the unstable manifold is expressed as a convergent power series of the form
\begin{equation}
\bm{X}_s = \varphi (\bm{X}_u) = \sum^\infty_{|n|=2} \bm{b}_n \bm{X}^n_u,
\end{equation}
where $n$ denotes a multi-index as usual.
The coefficient vectors $\bm{b}_n$ can be obtained by substituting it into Eq.(\ref{3-13}).
The system on the unstable manifold is given by
\begin{equation}
\left\{ \begin{array}{l}
\displaystyle \frac{d\bm{X}_u}{dt} = J_u\bm{X}_u + \bm{F}_u (\varphi (\bm{X}_u), \bm{X}_u, Z,\varepsilon ), \\[0.2cm]
\displaystyle \frac{dZ}{dt} = rZ +  \varepsilon F_{m+1}(\varphi (\bm{X}_u), \bm{X}_u,Z,\varepsilon ),  \\[0.2cm]
\displaystyle \frac{d\varepsilon }{dt} = s\varepsilon,
\end{array} \right.
\label{3-15}
\end{equation}

\textbf{Step 3.} Calculate the normal form of the first equation of (\ref{3-15}) up to degree $N$ to be determined;

Due to the normal form theory, there exists a polynomial transformation (near identity transformation)
\begin{equation}
\bm{X}_u \mapsto \bm{Y} = h_1(\bm{X}_u, Z, \varepsilon ), \quad h_1(0) = 0, Dh_1(0) = \mathrm{id},
\end{equation}
of degree $N$, which can be exactly calculated for a finite $N$, such that the first equation of Eq.(\ref{3-15}) takes the form (see Prop.A.1)
\begin{equation}
\left\{ \begin{array}{l}
\displaystyle \frac{d\bm{Y}}{dt} = J_u\bm{Y} + \bm{G}_1(\bm{Y},Z,\varepsilon ) + \bm{G}_2(\bm{Y},Z,\varepsilon ), \quad \bm{Y}\in \C^{m-k}, \\[0.2cm]
\displaystyle \frac{dZ}{dt} = rZ +  \varepsilon F_{m+1},  \\[0.2cm]
\displaystyle \frac{d\varepsilon }{dt} = s\varepsilon,
\end{array} \right.
\label{3-17}
\end{equation}
where $\bm{G}_1$ consists of resonance terms up to degree $N$, and $\bm{G}_2 \sim O(|| (\bm{Y},Z,\varepsilon ) ||^{N+1})$.
If we assume the condition (iii) of Prop.3.5, $\bm{G}_1 = 0$.
Nevertheless, we keep the term $\bm{G}_1$ to observe that what happen when a given system (\ref{3-1}) does not have the Painlev\'{e} property.
How to choose $N$ will be explained in Sec.4.2.

\textbf{Step 4.} Weighted blow-up;

Now the origin of the $m-k+2$ dimensional system (\ref{3-17}) on the unstable manifold is a singularity;
Laurent series solutions under consideration lie on the unstable manifold and they approach to the origin as $z\to z_0$.
In order to resolve the singularity, we introduce the weighted blow-up with the weight $(\lambda _{k+1}, \cdots , \lambda _m, r,s)$.
Roughly speaking, the weighted blow-up is a birational transformation $\pi : B \to \C^{m-k+2}$ whose exceptional divisor $\pi^{-1}(0)$
is the weighted projective space $\C P^{m-k+1} (\lambda _{k+1}, \cdots ,\lambda _m ,r ,s)$, and $B$ is a line bundle over 
$\C P^{m-k+1} (\lambda _{k+1}, \cdots ,\lambda _m ,r ,s)$.
Denote $\bm{Y} = (Y_{k+1}, \cdots ,Y_{m})$.
One of the local coordinates $(u_{k+1}, \cdots ,u_m, \zeta, w)$ of $B$ is defined by
\begin{equation}
\left(
\begin{array}{@{\,}c@{\,}}
Y_{k+1} \\
 \vdots \\
Y_{m} \\
Z \\
\varepsilon 
\end{array}
\right) = \left(
\begin{array}{@{\,}c@{\,}}
u_{k+1}w^{\lambda _{k+1}} \\
 \vdots \\
u_{m}w^{\lambda _{m}}\\
\zeta w^{r} \\
 w^{s}
\end{array}
\right).
\label{3-18}
\end{equation}
Hence,  $w$ denotes a coordinate on a fiber and $(u_{k+1}, \cdots ,u_m, \zeta)$
is the inhomogeneous coordinates of the chart $\C^{m-k+1}/\Z_s$ of $\C P^{m-k+1} (\lambda _{k+1}, \cdots ,\lambda _m ,r ,s)$.
In particular, the set $\{ w=0\} \subset \C P^{m-k+1} (\lambda _{k+1}, \cdots ,\lambda _m ,r ,s)$ is attached at infinity of the 
original chart $\C^{m+1} = \{ (x_1,\cdots ,x_m, z)\}$.

The coordinate transformation between the original chart and the new coordinates
$(\bm{X}_s, u_{k+1}, \cdots ,u_m, \zeta,w)$ is given by
\begin{equation}
\left\{ \begin{array}{l}
x_i = (c_ic_j^{-p_i/p_j} + h_2(\bm{X}_s, u_{k+1}, \cdots ,u_m, \zeta ,w))w^{-p_i}, \quad (i\neq j)  \\
x_j = w^{-p_j},  \\
z = \zeta,
\end{array} \right.
\label{3-19}
\end{equation}
where $h_2$ is a polynomial with $h_2(0) = 0$ that is obtained by a finite step.

Due to the orbifold structure of the exceptional divisor $\C P^{m-k+1} (\lambda _{k+1}, \cdots ,s)$, the $\Z_s$ action
\begin{eqnarray*}
(u_{k+1}, \cdots ,u_m, \zeta, w) \mapsto (\omega ^{\lambda _{k+1}}u_{k+1}, \cdots ,\omega ^{\lambda _{m}}u_{m}, \omega ^r\zeta, \omega ^{-1}w),
\quad \omega :=e^{2\pi i /s},
\end{eqnarray*}
acts on the space $\{(u_{k+1}, \cdots ,u_m, \zeta,w) \}$.
This is compatible with the $\Z_s$ action (\ref{2-19}) of the original chart; the one action induces the other through the transformation (\ref{3-19}).

Put $J_u = \mathrm{diag} (\lambda _{k+1}, \cdots ,\lambda _m ) $ and $\bm{G}_i = (G_{i,k+1}, \cdots ,\allowbreak  G_{i,m})$.
By the blow-up (\ref{3-18}), Eq.(\ref{3-17}) is transformed into the system
\begin{eqnarray*}
\left\{ \begin{array}{l}
\displaystyle \frac{du_i}{dt} = w^{-\lambda _i} (G_{1,n} + G_{2,n}), \\[0.2cm]
\displaystyle \frac{d\zeta}{dt} = w F_{m+1},  \\[0.2cm]
\displaystyle \frac{dw}{dt} = w.
\end{array} \right.
\end{eqnarray*}
Using $\zeta = z$ and deleting $t$, we obtain the system
\begin{equation}
\left\{ \begin{array}{l}
\displaystyle \frac{du_{i}}{dz} = w^{-1-\lambda _{i}}(G_{1,i} + G_{2,i})/F_{m+1}, \quad i=k+1,\cdots ,m,  \\[0.2cm]
\displaystyle \frac{dw}{dz} =1/F_{m+1}.   \\
\end{array} \right.
\label{3-20}
\end{equation}
Since $1/F_{m+1} = -(f_j + \varepsilon G_j)/p_j$ is holomorphic in $u_{k+1}, \cdots ,u_m, w$ and $z$, a singularity of the right hand side may
arise only from the factor $w^{-1-\lambda _{i}}$.
\\[0.2cm]
\textbf{Proposition \thedef.} 
If $N$ is sufficiently large, the function $w^{-1-\lambda _{i}}G_{2,i}$ is holomorphic in $u_{k+1},\cdots ,u_m, w, z$, while
$w^{-1-\lambda _{i}}G_{1,i}$ is of the form 
\begin{eqnarray*}
w^{-1-\lambda _{i}}G_{1,i} = w^{-1} \times \text{(polynomial of } u_{k+1},\cdots ,u_m, z).
\end{eqnarray*}
If the conditions of Prop.3.5 are satisfied, $G_{1,n} = 0$ and the right hand side of Eq.(\ref{3-20})
is holomorphic in $u_{k+1},\cdots ,u_m, w, z$.
Further, $1/F_{m+1} \neq 0$ when $w=0$.
Hence, there are no singularities of the foliation on the exceptional divisor $\{ w=0\}$.
\\[-0.2cm]

As a result, the singularity of the foliation at the point (\ref{3-5}) is resolved;
$m-k+1$-parameter family of integral curves that lie on the unstable manifold, all of which pass through the fixed point (\ref{3-5})
in $(X_1,\cdots ,X_m, Z, \varepsilon )$ coordinates, 
intersect with the $m-k+1$-dimensional exceptional divisor $\{w = 0\}$ at different points.
Further, if all K-exponents other than $-1$ are positive integers, then the right hand of (\ref{3-20}) is polynomial because 
a transcendental function may arise only from the expression of the unstable manifold $\bm{X}_s = \varphi (\bm{X}_u)$.
 \\[0.2cm]
\textbf{Proof.}  
Since
\begin{eqnarray*}
\bm{G}_2 (\bm{Y},Z,\varepsilon ) = \bm{G}_2 (u_{k+1}w^{\lambda _{k+1}},\cdots ,u_{m}w^{\lambda _m}, \zeta w^{r}, w^s)
\end{eqnarray*}
and $\bm{G}_2 \sim O(|| (\bm{Y},Z,\varepsilon ) ||^{N+1})$, $w^{-1-\lambda _{i}}G_{2,i}$ is holomorphic in $w$ if $N$ is sufficiently large.
On the other hand, since $G_{1,i}(\bm{Y},Z,\varepsilon )$ consists of resonance terms,
a monomial $\alpha := Y_{k+1}^{n_{k+1}}\cdots Y_{m}^{n_m}Z^{n_{m+1}}\varepsilon ^{n_{m+2}}$ included in $G_{1,i}(\bm{Y},Z,\varepsilon )$ satisfies
\begin{eqnarray*}
\lambda _{k+1} n_{k+1} + \cdots + \lambda _{m} n_{m} + r n_{m+1} + sn_{m+2} = \lambda _{i}.
\end{eqnarray*}
Hence, $w^{-1-\lambda _{i}}\alpha $ becomes of order $1/w$ after the blow-up.

Let us confirm the last statement.
When $w=0$, we have $Y_i = Z = \varepsilon =0$.
This implies $X_i = c_ic_j^{-p_i/p_j}$ when $w=0$.
Recall that $f_j$ in Eq.(\ref{3-3}) implies $f_j = f_j(X_1, \cdots ,1, \cdots ,X_m, Z)$.
Therefore, we obtain
\begin{eqnarray*}
1/F_{m+1}|_{w=0} &=& -\frac{1}{p_j}(f_j + \varepsilon G_j)|_{w=0} \\
&=&  -\frac{1}{p_j}f_j^A(c_1c_j^{-p_1/p_j}, \cdots ,1 ,\cdots , c_mc_j^{-p_m/p_j}) \\
&=& -\frac{1}{p_j} c_j^{-(1+p_j)/p_j} f_j^A(c_1, \cdots ,c_m) =c_j^{-1/p_j} \neq 0.\quad \Box
\end{eqnarray*}

\textbf{Step 5.} Divide (\ref{3-19}) and (\ref{3-20}) by the $\Z_{p_j}$-action;

If $p_j\neq 1$, (\ref{3-19}) is not a one-to-one transformation.
Recall that the group $\Z_{p_j}$
\begin{eqnarray*}
(X_i, Z, \varepsilon ) \mapsto (e^{2\pi i \cdot p_i/p_j}X_i, \,e^{2\pi i\cdot r/p_j}Z, \,e^{2\pi i\cdot s/p_j} \varepsilon ), \quad (i\neq j)
\end{eqnarray*}
acts on the inhomogeneous coordinates on the lift of the chart $\C^{m+1}/\Z_{p_j}$
and Eq.(\ref{3-3}) is invariant under the action due to the orbifold structure.
This action induces a $\Z_{p_j}$ action on the $(\bm{X}_s, u_{k+1}, \cdots ,w)$-coordinates and Eq.(\ref{3-20}) is invariant under the action;
\begin{eqnarray*}
\Z_{p_j} \curvearrowright \C^{m}_1 := \{ (\bm{X}_s, u_{k+1}, \cdots , u_m, w)\}.
\end{eqnarray*}
Obviously, the right hand sides of the transformation (\ref{2-25}) are invariant under the $\Z_{p_j}$ action.
This shows that the right hand sides of the transformation (\ref{3-19}) are rational invariants of the $\Z_{p_j}$ action.
Thus, if we divide $\C^{m}_1$ by the action, (\ref{3-19}) becomes a one-to-one rational transformation,
which can be explicitly given by rewriting the right hand sides of (\ref{3-19}) in terms of  polynomial invariants 
of the action $\Z_{p_j} \curvearrowright \C^{m+1}_1$, one of which should be $W := w^{p_j}$.
The original chart $M_0 := \{ (x_1, \cdots ,x_m)\} \simeq \C^{m}$ and the quotient space $M_1 :=\C^{m}_1/\Z_{p_j}$, which is a
nonsingular algebraic variety, are glued by the one-to-one rational transformation (\ref{3-19}) to give a nonsingular algebraic variety $M_{01}$.
Eq.(\ref{3-1}) together with (\ref{3-20}) gives a holomorphic equation on $M_{01}$ without singularities of the foliation. 

Suppose that a given system with the Painlev\'{e} property (in the sense that any solutions are meromorphic)
has $n$-types Laurent series solutions, all of whose leading terms are of the form $x_i \sim c_i(z-z_0)^{-p_i}$
; that is, there are $n$ roots of the equation $-p_ic_i = f_i^A(c_1, \cdots ,c_m)$. 
We perform Step 1 to Step 5 for all Laurent series to obtain the manifold $M_i \simeq \C^{m}/\Z_{p_j}$ and a holomorphic
differential equation on it as in Step 5.
Then, an algebraic variety $\mathcal{M}(z) := M_0 \cup M_1 \cup \cdots \cup M_n$ parameterized by $z \in \C$ 
gives the space of initial conditions for (\ref{3-1}).
Each solution defines a global holomorphic section of the fiber bundle $\{ (x,z) \, | \, x\in \mathcal{M}(z), z\in \C \}$
and there are no singularities of the foliation on the bundle.
See Chiba \cite{Chi2} for the detailed calculation for the first, second and fourth Painlev\'{e} equations, and Section 4 for the higher order
first Painlev\'{e} equation.


\section{The first Painlev\'{e} hierarchy}

Define the operator $\mathcal{L}_m$ by
\begin{eqnarray}
\frac{d}{dz}\mathcal{L}_{m+1}[x] = \left( \frac{d^3}{dz^3}-8x\frac{d}{dx}-4\frac{dx}{dz} \right) \mathcal{L}_m[x],
\quad \mathcal{L}_0[x] = 1,
\label{4-1}
\end{eqnarray}
where $x=x(z)$ is a function of $z\in \C$.
The $2m$-th order first Painlev\'{e} equation (the first Painlev\'{e} hierarchy) is defined to be $\mathcal{L}_m[x] = -4z$.
Indeed, it is easy to verify that there is a polynomial $P_m$ such that the equation is expressed as
\begin{eqnarray}
x^{(2m)} = P_m(x,x',\cdots ,x^{(2m-2)}) + z, \quad x^{(i)}:= \frac{d^ix}{dz^i}.
\label{4-2}
\end{eqnarray}
For example, we obtain
\begin{eqnarray*}
& & x'' = 6x^2 + z, \\
& &  x'''' = 20xx''+ 10(x')^2-40x^3+z,
\end{eqnarray*}
for $m=1,2$, respectively.
We rewrite (\ref{4-2}) as the $2m$-dimensional system 
\begin{eqnarray}
\text{($\text{P}_\text{I}$)}_m \left\{ \begin{array}{l}
x_1' = x_2  \\
\quad \vdots  \\
x_{2m-1}' = x_{2m} \\
x_{2m}' = P_m(x_1,\cdots ,x_{2m-1}) + z.
\end{array} \right.
\label{4-3}
\end{eqnarray}
This system satisfies the assumptions (A1) to (A3) with $g_i = 0$ and 
\begin{eqnarray*}
(p_1, p_2, \cdots ,p_{2m}, r, s) = (2,3,\cdots ,2m+1, 2m+2, 2m+3).
\end{eqnarray*}
Indeed, it is easy to prove by induction that the function $P_m$ satisfies
\begin{equation}
P_m (\lambda ^2 x_1,\lambda^3 x_2, \cdots , \lambda ^{2m}x_{2m-1} ) = \lambda ^{2m+2} P_m(x_1, x_2, \cdots ,x_{2m-1})
\label{4-4}
\end{equation}
for any $\lambda \in \C$.
Thus, the system (\ref{4-3}) induces a rational ODE on the weighted projective space $\C P^{2m+1}(2, \cdots ,2m+3)$.
In Shimomura\cite{Sim2}, it is proved that the first Painlev\'{e} hierarchy has the Painlev\'{e} property in the sense that
any solutions are meromorphic functions.
The leading coefficients of the Laurent series solutions are given by
\begin{eqnarray}
c_j(k) = (-1)^{j+1}j!\cdot b_0, \quad b_0:= \frac{1}{2}k(k+1),
\label{4-5}
\end{eqnarray}
for $k= 1, \cdots ,m$ \cite{Sim1}.
Hence, the system (\ref{4-3}) has $m$ Laurent series solutions of the form
\begin{eqnarray*}
x_i(z) = c_i(k) (z-z_0)^{-p_i} + \sum^\infty_{n=1} a_{i,n}(k) (z-z_0)^{-p_i+n}. 
\end{eqnarray*}


\subsection{The K-exponents of the first Painlev\'{e} hierarchy}

For each Laurent series solution with (\ref{4-5}), the K-exponents are defined, which are given as follows.
\\[0.2cm]
\textbf{Theorem \thedef.} The K-exponents of the system (\ref{4-3}) associated with the Laurent series solution with (\ref{4-5})
are given by the following $2m$ integers;
\begin{eqnarray*}
\lambda &=&\quad \,\, 2, \quad \quad 4, \quad\cdots ,2m-2k, \quad (m-k) \\
& & 2k+3,\, 2k+5,\cdots ,\,2m+1, \quad (m-k) \\
& & 2m+4,\, 2m+6,\cdots ,2m+2k+2, \quad (k) \\
& & \quad -1, \quad\,\,\, -3,\quad \cdots ,-(2k-1),\quad\quad \,\,(k) 
\end{eqnarray*}
Thus, the Laurent series solution includes $2m-k+1$ free parameters (including $z_0$).
In particular, the Laurent series for the case $k=1$ includes $2m$ free parameters that represents a general solution.
\\[-0.2cm]

In order to prove the theorem, we need a Hamiltonian form of the system.
By putting $x \mapsto \lambda ^2x$ and $z\mapsto \lambda ^{-1}z$ with $\lambda ^{-2m-3} =4^m$, Eq.(\ref{4-2}) is transformed to
the equation $x^{(2m)} = P_m(x,\cdots ,x^{(2m-2)}) + 4^m z$ due to (\ref{4-4}).
Further, by putting
\begin{eqnarray}
\left\{ \begin{array}{l}
\displaystyle u_j = 4^{1-j} (x_{2j-1} - P_{j-1} (x_1 ,\cdots ,x_{2j-3})), \\
\displaystyle v_j = \frac{4^{1-j}}{2}\left( x_{2j} - \sum^{2j-3}_{i=1} \frac{\partial P_{j-1}}{\partial x_i} (x_1 ,\cdots ,x_{2j-3}) x_{i+1}\right), \\
\end{array} \right.
\label{4-6}
\end{eqnarray}
$(u_j, v_j)$ satisfies the system
\begin{equation}
\left\{ \begin{array}{l}
u_j' = 2v_j,   \\
v_j' = 2u_{j+1} + 2u_1 u_j + 2w_j, \quad (j=1,\cdots ,m),  \\
\end{array} \right.
\label{4-7}
\end{equation}
where $u_{m+1} = 0$ and $w_j$ is determined by the recursive relation
\begin{eqnarray*}
w_j = \frac{1}{2}\sum^j_{k=1}u_ku_{j+1-k} + \sum^{j-1}_{k=1} u_kw_{j-k} - \frac{1}{2} \sum^{j-1}_{k=1} v_kv_{j-k} + \delta _{jm}z.
\end{eqnarray*}
The system (\ref{4-7}) is introduced by Shimomura \cite{Sim2} to prove the Painlev\'{e} property.
If we define the weighted degree of $x_j$ by $\mathrm{deg} (x_j) = j+1$, then Eqs.(\ref{4-4}) and (\ref{4-6}) provide
$\mathrm{deg} (u_j) = 2j$ and $\mathrm{deg} (v_j) = 2j+1$.
This implies that the transformation
\begin{eqnarray*}
(x_1 , \cdots ,x_{2m}, z) \mapsto (u_1, v_1 ,\cdots ,u_m, v_m, z)
\end{eqnarray*}
defined by (\ref{4-6}) is an automorphism on $\C P^{2m+1}(2, \cdots ,2m+3)$.
In particular, K-exponents of Eq.(\ref{4-3}) are the same as those of (\ref{4-7}) due to Thm.2.5 or Thm.3.4.

According to Takei \cite{Tak}, we further change coordinates by
\begin{equation}
u_j = (-1)^{j-1} Q_j, \quad v_j = 2(-1)^{m-j} (P_{m-j+1} + P_{m-j+2}Q_1 + \cdots +P_mQ_{j-1}).
\end{equation}
Then, $(P_j, Q_j)$ satisfies the Hamiltonian system
\begin{equation}
\frac{dP_j}{dz} = -\frac{\partial H_m}{\partial Q_j} ,\quad \frac{dQ_j}{dz} = \frac{\partial H_m}{\partial P_j}, \quad (j=1,\cdots ,m),
\label{4-9}
\end{equation}
where $H_m$ is a polynomial Hamiltonian function.
Since the weighted degrees are given by $\mathrm{deg} (P_j) = 2m+3-2j$ and $\mathrm{deg} (Q_j) = 2j$, the transformation
\begin{eqnarray*}
(u_1, v_1, \cdots ,u_m, v_m,z) \mapsto (P_1, Q_1, \cdots ,P_m, Q_m, z)
\end{eqnarray*}
is an isomorphism from $\C P^{2m+1}(2, \cdots ,2m+3)$ to
\begin{eqnarray*}
\C P^{2m+1}(2m+1,2,\cdots ,2m+3-2j,2j, \cdots ,3,2m,2m+2,2m+3).
\end{eqnarray*}
In particular, the K-exponents do not change.
It is easy to verify the equality
\begin{eqnarray}
H_m(\cdots ,\lambda ^{2m+3-2j}P_j, \lambda ^{2j}Q_j  ,\cdots , \lambda ^{2m+2}z)
 = \lambda ^{2m+4}H_m(\cdots ,P_j, Q_j, \cdots , z).
\label{4-10}
\end{eqnarray}
Thus Lemma 2.4 provides
\\[0.2cm]
\textbf{Lemma \thedef.} If $\lambda $ is a K-exponent of the system (\ref{4-3}), so is $\mu = 2m+3-\lambda $.
\\[-0.2cm]

Because of this lemma, the existence of K-exponents in the fourth line in Thm.4.1 immediately follows from that of the third line.
\\[0.2cm]
\textbf{Proof of Thm.4.1.} The K-matrix of the system (\ref{4-3}) is given by
\begin{eqnarray*}
K = \left(
\begin{array}{@{\,}cccccccc@{\,}}
2 & 1 & \cdots  & 0 & 0 & \cdots  & 0 & 0 \\
\vdots &\vdots & &\vdots &\vdots & & \vdots & \vdots \\
0 & 0 & \cdots  & j+1 &1 &\cdots  & 0 & 0 \\
\vdots &\vdots & &\vdots &\vdots & & \vdots & \vdots \\
0& 0 &\cdots  &0 &0 & \cdots & 2m & 1 \\
\frac{\partial P_m}{\partial x_1}&\frac{\partial P_m}{\partial x_2} &\cdots &
\frac{\partial P_m}{\partial x_j} &\frac{\partial P_m}{\partial x_{j+1}} & \cdots  & \frac{\partial P_m}{\partial x_{2m-1}} & 2m+1 \\
\end{array}
\right),
\end{eqnarray*}
where $\partial P_m/\partial x_j$ is estimated at the point
\begin{eqnarray*}
\bm{c}(k) = (c_1(k), \cdots , c_{2m-1}(k)).
\end{eqnarray*}
The eigen-equation is given by
\begin{eqnarray*}
\det (\lambda -K) = (\lambda -2) \cdots (\lambda -2m-1) - \sum^{2m-1}_{i=1}\frac{\partial P_m}{\partial x_i}(\bm{c}(k))(\lambda -2)\cdots (\lambda -i) = 0.
\end{eqnarray*}
By the definition, $P_m$ satisfies
\begin{eqnarray*}
\mathcal{L}_{m+1}[x] = -4 \left( x^{(2m)} - P_m(x, \cdots ,x^{(2m-2)}) \right) ,\quad \mathcal{L}_0[x] = 1.
\end{eqnarray*}
Putting $x = \varphi _0 + \delta \varphi _1$ yields
\begin{eqnarray*}
\mathcal{L}_{m+1}[ \varphi _0 + \delta \varphi _1]
    &=& -4 \left( \varphi _0^{(2m)} - P_m(\varphi _0, \cdots ,\varphi _0^{(2m-2)}) \right) \\
& & \quad -4 \delta \left( \varphi _1^{(2m)} - \sum^{2m-1}_{i=1} \frac{\partial P_m}{\partial x_i}
    (\varphi _0, \cdots ,\varphi _0^{(2m-2)}) \varphi _1^{(i-1)}\right) + O(\delta ^2). 
\end{eqnarray*}
If we put $\varphi _0(z) = b_0(z+1)^{-2}$ with $b_0 = k(k+1)/2$, then the first term in the right hand side vanishes because of Eq.(\ref{4-4}).
Since
\begin{eqnarray*}
\varphi _0^{(j)}(0) = (-1)^j (j+1)! b_0 = c_{j+1}(k),
\end{eqnarray*}
we obtain
\begin{eqnarray*}
\mathcal{L}_{m+1}[ \varphi _0 + \delta \varphi _1](0)
   = -4 \delta \left( \varphi _1^{(2m)}(0) - \sum^{2m-1}_{i=1} \frac{\partial P_m}{\partial x_i}(\bm{c}(k)) \varphi _1^{(i-1)}(0) \right) + O(\delta ^2). 
\end{eqnarray*}
Further, putting $\varphi _1(z) = b_0 (z+1)^{\lambda -2}$ provides
\begin{eqnarray*}
\mathcal{L}_{m+1}[ \varphi _0 + \delta \varphi _1](0) = -4 \delta b_0 \cdot \det (\lambda -K) + O(\delta ^2). 
\end{eqnarray*}
Therefore, if we set
\begin{equation}
\mathcal{L}_{j}[ \varphi _0 + \delta \varphi _1](z) = f_j(z) + \delta g_j(z) + O(\delta ^2), 
\quad f_0 = 1,\, g_0 = 0,
\label{4-11}
\end{equation}
$\det (\lambda -K)=0$ is equivalent to $g_{m+1}(0) = 0$.

Let us derive difference equations for $f_j$ and $g_j$.
Substituting (\ref{4-11}) into the definition (\ref{4-1}) of $\mathcal{L}_{j+1}$, we obtain
\begin{equation}
\left\{ \begin{array}{l}
f_{j+1}' = f_j''' - 8\varphi _0f_j'-4\varphi _0'f_j,  \\
g_{j+1}' = g_j''' - 8\varphi _0g_j'-4\varphi _0'g_j - 8\varphi _1f_j' - 4\varphi _1' f_j.  \\
\end{array} \right.
\label{4-12}
\end{equation}
If we set $f_j = A_j(z+1)^{-2j}$, the first equation yields
\begin{eqnarray*}
(2j+2)A_{j+1} &=& 2j(2j+1)(2j+2)A_{j} - 16 jb_0A_j - 8b_0A_j \\
&=& 8(2j+1)\left( \frac{1}{2}j(j+1) - b_0 \right) A_j.
\end{eqnarray*}
Thus, we have
\begin{equation}
A_{j+1} = -4b_0 \prod^j_{l=1}\frac{4(2l+1)}{l+1}\left( \frac{1}{2}l(l+1)-b_0\right), \quad b_0 = \frac{1}{2}k(k+1).
\end{equation}
This is further rearranged as
\begin{equation}
A_j = (-1)^j 2^j \frac{(2j-1)!! \cdot (k+j)!}{j! \cdot (k-j)!},\quad (j=1,\cdots ,k),
\label{4-14}
\end{equation}
and $A_j = 0$ for $j\geq k+1$.

Next, by putting $g_j = B_j (z+1)^{\lambda -2j}$, the second equation of (\ref{4-12}) gives
\begin{eqnarray*}
B_{j+1}=\frac{\left( \lambda -(2j+2k+2)\right) \left( \lambda -(2j+1)\right) \left( \lambda -(2j-2k) \right)}{\lambda -(2j+2)} B_j
 - \frac{4 \left( \lambda -(4j+2)\right)}{\lambda -(2j+2)} b_0 A_j.
\end{eqnarray*}
Since $g_{m+1}(0) = 0$ if and only if $B_{m+1} = 0$, roots of $B_{m+1}(\lambda ) = 0$ give the K-exponents.
Since $A_j = 0$ for $j\geq k+1$, we obtain
\begin{eqnarray}
& & B_{m+1} =\frac{\left( \lambda -(2m+2k+2)\right) \left( \lambda -(2m+1)\right) \left( \lambda -(2m-2k) \right)}{\lambda -(2m+2)} B_m \nonumber \\
&=& \frac{\prod^{m}_{j=k+1} \left( \lambda -(2j+2k+2)\right)  \cdot \left( \lambda -(2j+1)\right)  \cdot\left( \lambda -(2j-2k) \right)}
{(\lambda -(2m+2)) \cdot (\lambda -2m)\cdots (\lambda -(2k+4))} B_{k+1}.
\label{4-15}
\end{eqnarray}
Now we need two lemmas.
\\[0.2cm]
\textbf{Lemma \thedef.} For any $j\geq k$, $B_{j+1}(\lambda )$ is a polynomial in $\lambda $ of degree $2 j$. 
\\[0.2cm]
\textbf{Lemma \thedef.} The equation $B_{k+1}(\lambda ) = 0$ has $k$ roots given by
$\lambda = 2k+4,\,\, 2k+6, \cdots , 4k+2$.
In particular, there is a polynomial $C_{k+1}(\lambda )$ of degree $k$ such that 
\begin{equation}
B_{k+1} = (\lambda - (2k+4)) \cdots (\lambda -(4k+2)) C_{k+1}(\lambda ).
\end{equation} 
Lemma 4.3 is trivial because $B_{j+1}(\lambda )=0$ is equivalent to the eigen-equation $\det (\lambda - K) = 0$
for the $2j$ dimensional problem.
If Lemma 4.4 is true, all factors in the denominator of (\ref{4-15}) cancel and we obtain
\begin{eqnarray*}
B_{m+1} =\!\! \prod^{m}_{j=m+1-k}\!\!\!\! (\lambda -(2j+2k+2)) \!\!\prod^{m}_{j=k+1} \!\!(\lambda -(2j+1)) \!\!\prod^{m}_{j=k+1} \!\!(\lambda -(2j-2k)) \cdot C_{k+1}(\lambda ).
\end{eqnarray*}
In particular, we obtained the first three lines in Thm.4.1.
This completes the proof because of Lemma 4.2.
\\[0.2cm]
\textbf{Proof of Lemma 4.5.} 
Put 
\begin{eqnarray*}
P_j = \frac{\left( \lambda -(2j+2k+2)\right) \left( \lambda -(2j+1)\right) \left( \lambda -(2j-2k) \right)}{\lambda -(2j+2)},
\quad Q_j = -\frac{2 \left( \lambda -(4j+2)\right)}{\lambda -(2j+2)}.
\end{eqnarray*} 
Then we have
\begin{eqnarray*}
B_{k+1} &=& P_k B_k + Q_k k(k+1) A_k \\
&=& P_k(P_{k-1}B_{k-1}+Q_{k-1} k(k+1) A_{k-1}) + Q_k k(k+1) A_k\\
&\vdots &  \\
&=&P_kP_{k-1}\cdots P_1B_1+ [Q_kA_k + P_kQ_{k-1}A_{k-1} +\cdots +P_k\cdots P_2Q_1A_1] k(k+1).
\end{eqnarray*}
Substituting $P_j, Q_j$ and (\ref{4-14}), we obtain
\begin{eqnarray*}
\frac{B_{k+1}(\lambda )}{k(k+1)} &=& \sum^k_{l=0}\frac{\lambda -(4k-4l+2)}{\lambda -(2k-2l+2)}(-1)^{k-l+1} 2^{k-l+1}
\frac{(2k-2l-1)!! \cdot (2k-l)!}{(k-l)! \cdot l!} \times \\
& & \qquad \prod^k_{j=k-l+1} \frac{\left( \lambda -(2j+2k+2)\right) \left( \lambda -(2j+1)\right) \left( \lambda -(2j-2k) \right)}{\lambda -(2j+2)}.
\end{eqnarray*}
Now we show that $\lambda = 4k+2-2n$ is a root of $B_{k+1}(\lambda )=0$ for $n=0,1,\cdots ,k-1$.
Substituting this value gives
\begin{eqnarray*}
\frac{B_{k+1}(\lambda )}{k(k+1)} &=& \sum^k_{l=0}\frac{2l-n}{k+l-n}(-1)^{k-l+1} 2^{k-l+1}
\frac{(2k-2l-1)!! \cdot (2k-l)!}{(k-l)! \cdot l!} \times \\
& & \qquad \prod^k_{j=k-l+1} \frac{(-2)\left( j+n-k \right) \left( 4k-2n-2j+1 \right) \left( 3k-n-j+1 \right)}{2k-n-j}.
\end{eqnarray*}
Since the factor $j+n-k$ becomes zero when $j=k-n$, which is possible only for $l=n+1, \cdots ,k$, 
\begin{eqnarray*}
\frac{B_{k+1}(\lambda )}{k(k+1)} &=& (-1)^k 2^{k+1} \sum^n_{l=0}\frac{2l-n}{k+l-n}
\frac{(2k-2l-1)!! \cdot (2k-l)!}{(k-l)! \cdot l!} \times \\
& & \qquad \prod^k_{j=k-l+1} \frac{\left( j+n-k \right) \left( 4k-2n-2j+1 \right) \left( 3k-n-j+1 \right)}{2k-n-j}.
\end{eqnarray*}
Define
\begin{eqnarray*}
F(l)&:=& \frac{2l-n}{k+l-n}
\frac{(2k-2l-1)!! \cdot (2k-l)!}{(k-l)! \cdot l!} \times \\
& & \qquad
\prod^k_{j=k-l+1} \frac{\left( j+n-k \right) \left( 4k-2n-2j+1 \right) \left( 3k-n-j+1 \right)}{2k-n-j}.
\end{eqnarray*}
Then, it is straightforward to prove that $F(l) = -F(n-l)$ and $F(n/2) = 0$ when $n$ is an even number.
This proves $B_{k+1}(\lambda ) = 0$ for $\lambda = 4k+2-2n$ with $n=0,1,\cdots ,k-1$. $\Box$


\subsection{The space of initial conditions for the fourth order equation}

The fourth order first Painlev\'{e} equation is given by
\begin{eqnarray}
\text{($\text{P}_\text{I}$)}_2 \left\{ \begin{array}{l}
x_1' = x_2  \\
x_2' = x_3 \\
x_3' = x_4 \\
x_4' =20x_1 x_3 + 10x_2^2 - 40x_1^3 + z.
\end{array} \right.
\label{4-17}
\end{eqnarray}
In this section, we demonstrate how to construct the space of the initial conditions of this system.
The system satisfies the assumptions (A1) to (A3) with the weight
\begin{equation}
(p_1, p_2, p_3, p_4, r,s) = (2,3,4,5,6,7).
\end{equation}
Thus, we give the system on the local chart $\C^5/\Z_7$ of the space $\C P^5 (2,3, \cdots ,7)$.
The system has the two families of Laurent series solutions 
\begin{eqnarray}
\mathrm{(I)}& &  x_j (z) \sim (-1)^{j+1} j! \,(z-z_0)^{-p_j}, \label{4-19} \\
\mathrm{(II)}& & x_j (z) \sim 3(-1)^{j+1} j! \,(z-z_0)^{-p_j},\label{4-20}
\end{eqnarray}
whose K-exponents are given by
\begin{eqnarray*}
\mathrm{(I)}& &  \lambda = -1, 2, 5, 8, \\
\mathrm{(II)}& & \lambda = -1, -3, 8, 10,
\end{eqnarray*}
respectively.
In particular, the first one represents a general solution.
We perform the resolution of singularity for each Laurent series.
\\

\textbf{(I)} Let us consider the resolution of the Laurent series (I).

\textbf{Step 1.} The coordinate transformation between the original coordinates and the inhomogeneous coordinates on $\C^5/ \Z_2$ are given by
\begin{equation}
x_1 = \varepsilon ^{-2/7},\,\, x_2 = X_2 \varepsilon ^{-3/7},\, x_3 = X_3 \varepsilon ^{-4/7},\, x_4 = X_4 \varepsilon ^{-5/7},\, z = Z \varepsilon ^{-6/7}.
\end{equation}
We express the system (\ref{4-17}) in the new coordinates as a polynomial vector field of the form (\ref{3-3c})
\begin{equation}
\left\{ \begin{array}{ll}
\displaystyle dX_2/dt = 3X_2^2 - 2X_3,  \\
\displaystyle dX_3/dt = 4X_3X_2 - 2X_4,  \\
\displaystyle dX_4/dt = 5X_4X_2 - (40 X_3 +20 X_2^2-80+2z),  \\
\displaystyle dZ/dt = 6ZX_2 - 2\varepsilon,  \\
\displaystyle d\varepsilon /dt = 7\varepsilon X_2.   \\
\end{array} \right.
\label{4-22}
\end{equation}
This system has two fixed points
\begin{eqnarray*}
(X_2, X_3, X_4, Z, \varepsilon ) = (2,\, 6,\, 24,\, 0,0), \quad (2/\sqrt{3},\,\, 2,\,\, 8/\sqrt{3},\,\, 0,0),
\end{eqnarray*}
which correspond to the Laurent series (I) and (II), respectively.
The eigenvalues of the Jacobi matrix at these fixed points are 
\begin{eqnarray*}
& & \lambda = 4,\,\, 10,\,\, 16,\,\, 12,\,\, 14, \\
& & \lambda = -2\sqrt{3},\,\,  16/\sqrt{3},\,\,  20/\sqrt{3},\,\,  12/ \sqrt{3},\,\,  14/\sqrt{3}, 
\end{eqnarray*}
respectively.
Since we have used the polynomial form (\ref{3-3c}) instead of (\ref{3-3}),
they differ from the K-exponents, $r$ and $s$ by a constant factor 
(multiplied by $1/2$ and $\sqrt{3}/2$, they become $2,5,8,6,7$ for the first one and $-3,8,10,6,7$ for the second one, respectively,
which coincide with the K-exponents and $r,s$).

For the resolution of singularity of the first fixed point (I), put 
\begin{equation}
(\hat{X}_2, \hat{X}_3, \hat{X}_4 ) = (X_2 - 2, X_3 -6, X_4 - 24)
\end{equation}
to obtain the system of the form (\ref{3-12}).

\textbf{Step 2.} We introduce the linear transformation
\begin{equation}
\left\{ \begin{array}{l}
\hat{X}_2 = v_1, \\
\hat{X}_3 = 4 v_1 + v_2, \\
\hat{X}_4 = 20 v_1 + 3 v_2 + v_3 + Z/2 - \varepsilon /2.
\end{array} \right.
\end{equation}
Then, we obtain the system of the form
\begin{equation}
\left\{ \begin{array}{l}
\displaystyle dv_1/dt = 4v_1 + F_1(v_1, v_2, v_3,Z, \varepsilon )  \\
\displaystyle dv_2/dt = 10v_2+ F_2(v_1, v_2, v_3,Z, \varepsilon )  \\
\displaystyle dv_3/dt = 16v_3 + F_3(v_1, v_2, v_3,Z, \varepsilon )  \\
\displaystyle dZ/dt = 12Z + F_4 (v_1, v_2, v_3,Z, \varepsilon ) \\
\displaystyle d\varepsilon /dt = 14\varepsilon + F_5(v_1, v_2, v_3,Z, \varepsilon ),  \\
\end{array} \right.
\label{4-25}
\end{equation}
where $F_1, \cdots , F_5$ are defined by
\begin{eqnarray*}
\left\{ \begin{array}{l}
F_1 = -2v_2 + 3v_1^2,  \\
F_2 = -2v_3 - Z + \varepsilon + 4v_1^2 + 4 v_1v_2,   \\
F_3 = 8 v_1^2 + 3v_1v_2 + 5 v_1v_3 - v_1Z/2 + v_1 \varepsilon , \\
F_4 = -2\varepsilon + 6v_1Z, \\
F_5 = 7v_1 \varepsilon .
\end{array} \right.
\end{eqnarray*}
Note that we need not diagonalize the linear part;
$-2v_2$ in $F_1$, $-2v_3 - Z + \varepsilon$ in $F_2$ and $-2\varepsilon$ in $F_4$ do not yield a singularity after the blow-up
(see the next step). 

\textbf{Step 3.} We calculate the normal form of the system (\ref{4-25}) to delete several monomials included in $F_1, F_2, F_3$.
We define weighted degrees to be
\begin{equation}
\mathrm{deg} (v_1) = 2, \,\,\mathrm{deg} (v_2) = 5, \,\,\mathrm{deg} (v_3) = 8, \,\,\mathrm{deg} (Z) = 6, \,\,\mathrm{deg} (\varepsilon ) = 7,
\end{equation}
which are the same as the weights of the weighted blow-up done in Step 4.
From the argument of Step 4 in Sec.3.3, it turns out that if a monomial $\alpha $ included in $F_i$
satisfies $\mathrm{deg} (\alpha ) < \mathrm{deg}(v_i)+1$, then the monomial yields a factor $1/w^n$ for some $n\geq 1$
in the right hand side of the system after the blow-up.
Hence, we have to remove such monomials by the normal form theory.
Monomials which may yield the factor $1/w^n$ are $v_1^2$ in $F_2$ and $v_1^2, v_1v_2, v_1^3, v_1^4$ in $F_3$
(although $v_1^3$ and $v_1^4$ are not included in $F_3$, they may appear after removing $v_1^2$ and $v_1v_2$).
To remove them, we set
\begin{equation}
\left\{ \begin{array}{l}
y_1 = v_1, \\
y_2 = v_2 + a_1 v_1^2,  \\
y_3 = v_3 + a_2 v_1^2 + a_3 v_1v_2 + a_4 v_1^3 + a_5 v_1^4.  \\
\end{array} \right.
\end{equation}
We can verify by a straightforward calculation that if we put
\begin{eqnarray*}
a_1 = 3,\,\, a_2 = 1,\,\, a_3 = -1/2,\,\, a_4 = -1/2,\,\, a_5 = 0,
\end{eqnarray*}
then the system (\ref{4-25}) is brought into
\begin{equation}
\left\{ \begin{array}{l}
\displaystyle dy_1/dt = 4y_1 - 2y_2 + 9 y_1^2,  \\
\displaystyle dy_2/dt = 10y_2 - 2 y_3 - Z + \varepsilon - 9y_1y_2 + 44 y_1^3,  \\
\displaystyle dy_3/dt = 16y_3 + 6 y_1 y_3 + y_1 \varepsilon /2+ y_2^2 - 7 y_1^2 y_2 /2  \\
\displaystyle dZ/dt = 12Z - 2\varepsilon + 6y_1Z \\
\displaystyle d\varepsilon /dt = 14\varepsilon + 7 y_1 \varepsilon. \\
\end{array} \right.
\label{4-28}
\end{equation}
This system does not include a monomial $\alpha $ satisfying $\mathrm{deg} (\alpha ) < \mathrm{deg}(v_i)+1$ except for
the diagonal part $(4y_1, 10y_2, 16 y_3, 12 Z, 14 \varepsilon )$, which will be removed by the blow-up below.

\textbf{Step 4.} We employ the weighted blow-up by
\begin{equation}
\left(
\begin{array}{@{\,}c@{\,}}
y_1 \\
y_2 \\
y_3 \\
Z \\
\varepsilon 
\end{array}
\right) = \left(
\begin{array}{@{\,}c@{\,}}
u_{1}w^{2} \\
u_{2}w^{5} \\
u_{3}w^{8} \\
\zeta w^{6} \\
 w^{7}
\end{array}
\right).
\label{4-29}
\end{equation}
Then, we obtain the polynomial system as desired;
\begin{equation}
\left\{ \begin{array}{l}
\displaystyle \frac{du_1}{dz} = u_2 w^2 - \frac{7}{2}u_1^2 w,  \\[0.3cm]
\displaystyle \frac{du_2}{dz} = -22u_1^3 + 7u_1u_2w + u_3 w^2 - \frac{1}{2}w + \frac{1}{2}z,  \\[0.3cm]
\displaystyle \frac{du_3}{dz} = u_1u_3w - \frac{1}{4}u_1 + \frac{7}{4}u_1^2u_2 - \frac{1}{2}u_2^2 w,  \\[0.3cm]
\displaystyle \frac{dw}{dz} =-1 - \frac{1}{2}u_1w^2.  \\[0.3cm]
\end{array} \right.
\end{equation}
The coordinate transformation between the original coordinates and $(u_1,u_2,u_3,w, \zeta)$ is given by
\begin{equation}
\left\{ \begin{array}{l}
\displaystyle x_1 = w^{-2}  \\
\displaystyle x_2 = \left( 2 + u_1w^2  \right)w^{-3},  \\
\displaystyle x_3 = \left( 6 + 4u_1w^2 - 3u_1^2w^4 + u_2w^5 \right)w^{-4}  \\
\displaystyle x_4 = \bigl( 24+20u_1w^2 - 10u_1^2 w^4 + 3u_2 w^5 - u_1^3 w^6 \\
\displaystyle \quad \quad \quad \quad \quad
 - \frac{1}{2}w^7 + \frac{1}{2}u_1u_2 w^7 + u_3 w^8 + \frac{1}{2}w^6z \bigr) w^{-5}, \\
\displaystyle z = \zeta.
\end{array} \right.
\label{4-30}
\end{equation}

\textbf{Step 5.} Due to the orbifold structure of $\C P^5 (2,3,4,5,6,7)$, the $\Z_2$ action
\begin{equation}
(X_2, X_3, X_4, Z, \varepsilon ) \mapsto (-X_2, X_3, -X_4, Z, -\varepsilon )
\label{Z2}
\end{equation}
acts on the coordinates $(X_2, X_3, X_4, Z, \varepsilon )$.
This induces the $\Z_2$ action on the $(u_1, u_2, u_3, z, w)$ coordinates given by
\begin{eqnarray}
 (u_1, u_2, u_3, z, w) &\mapsto & (-u_1-\frac{4}{w^2}, -u_2- \frac{32u_1}{w^3} - \frac{64}{w^5},\nonumber \\
& &       -u_3 - \frac{z}{w^2}- \frac{4 u_2}{w^3}+ \frac{24u_1^2}{w^4} + \frac{32u_1}{w^6}+ \frac{64}{w^8},z , -w). 
\end{eqnarray}
If we divide $\C_1^4 = \{ (u_1, u_2, u_3, w)\}$ by this action, (\ref{4-30}) becomes a one-to-one rational transformation
which defines a smooth algebraic variety $\C_0^4 \cup \C^4_1 /\Z_2$,
where $\C^4_0 = \{ (x_1,x_2,x_3,x_4)\}$ is the original chart.
\\

\textbf{(II)} Next, we consider the resolution of the Laurent series (II).

\textbf{Step 1.}
For the resolution of singularity of the second fixed point (II) of the system (\ref{4-22}), put 
\begin{equation}
(\hat{X}_2, \hat{X}_3, \hat{X}_4 ) = (X_2 - 2/\sqrt{3}, X_3 -2, X_4 - 8/\sqrt{3})
\end{equation}
to obtain the system of the form (\ref{3-12}).

\textbf{Step 2.} We introduce the linear transformation
\begin{equation}
\left\{ \begin{array}{l}
\hat{X}_2 = v_1,   \\
\hat{X}_3 = 3\sqrt{3} v_1 +  v_2 + \frac{3}{8}Z - \frac{5\sqrt{3}}{8}\varepsilon, \\
\hat{X}_4 = 25 v_1 + \frac{5\sqrt{3}}{3} v_2 + v_3 + \frac{7\sqrt{3}}{8}Z - \frac{27}{8}\varepsilon .
\end{array} \right.
\end{equation}
Then, we obtain the system of the form
\begin{eqnarray*}
\left\{ \begin{array}{l}
 dv_1/dt = -2\sqrt{3}v_1 -2v_2 - \frac{3}{4}Z + \frac{5\sqrt{3}}{4}\varepsilon  + F_1(v_1, v_2, v_3,Z, \varepsilon )  \\[0.2cm]
 dv_2/dt = \frac{16}{\sqrt{3}}v_2 -2 v_3 + F_2(v_1, v_2, v_3,Z, \varepsilon )  \\[0.2cm]
dv_3/dt = \frac{20}{\sqrt{3}}v_3 + F_3(v_1, v_2, v_3,Z, \varepsilon )  \\[0.2cm]
 dZ/dt =\frac{12}{\sqrt{3}}Z - 2\varepsilon + F_4 (v_1, v_2, v_3,Z, \varepsilon ) \\[0.2cm]
d\varepsilon /dt = \frac{14}{\sqrt{3}}\varepsilon + F_5(v_1, v_2, v_3,Z, \varepsilon ),  \\
\end{array} \right.
\end{eqnarray*}
where $F_1, \cdots , F_5$ are nonlinear terms.
The unstable manifold is a $(v_2, v_3, Z, \varepsilon )$-space.
We denote the unstable manifold by
\begin{equation}
v_1 = \varphi (v_2, v_3, Z, \varepsilon ),
\end{equation}
with a convergent power series $\varphi $ which does not include a constant term.
The system on the unstable manifold is given by
\begin{equation}
\left\{ \begin{array}{l}
dv_2/dt = \frac{16}{\sqrt{3}}v_2 -2 v_3 + F_2( \varphi (v_2, v_3, Z, \varepsilon ), v_2, v_3,Z, \varepsilon )  \\[0.2cm]
dv_3/dt = \frac{20}{\sqrt{3}}v_3 + F_3( \varphi (v_2, v_3, Z, \varepsilon ), v_2, v_3,Z, \varepsilon )  \\[0.2cm]
 dZ/dt = \frac{12}{\sqrt{3}}Z - 2\varepsilon + F_4 ( \varphi (v_2, v_3, Z, \varepsilon ), v_2, v_3,Z, \varepsilon ) \\[0.2cm]
d\varepsilon /dt = \frac{14}{\sqrt{3}}\varepsilon + F_5( \varphi (v_2, v_3, Z, \varepsilon ), v_2, v_3,Z, \varepsilon ).  \\
\end{array} \right.
\label{4-25b}
\end{equation}

\textbf{Step 3.}
We define weighted degrees by
\begin{equation}
\mathrm{deg} (v_2) = 8, \,\,\mathrm{deg} (v_3) = 10, \,\,\mathrm{deg} (Z) = 6, \,\,\mathrm{deg} (\varepsilon ) = 7,
\end{equation}
which are the same as the weights of the weighted blow-up done in Step 4.
As before, if a monomial $\alpha $ included in $F_i\, (i=2,3)$
satisfies $\mathrm{deg} (\alpha ) < \mathrm{deg}(v_i)+1$, then the monomial yields a factor $1/w^n$
in the right hand side of the system after the blow-up.
Since $F_i$ is nonlinear, $F_i( \varphi (v_2, v_3, Z, \varepsilon ), v_2, v_3,Z, \varepsilon )$
does not include such monomials (the possible least degree among nonlinear monomials is $\mathrm{deg} (Z^2) = 12$, which is larger than $\mathrm{deg} (v_2)+1$
and $\mathrm{deg} (v_3)+1$).
Hence, we need not calculate the normal form.

\textbf{Step 4.} We employ the weighted blow-up by
\begin{equation}
\left(
\begin{array}{@{\,}c@{\,}}
v_2 \\
v_3 \\
Z \\
\varepsilon 
\end{array}
\right) = \left(
\begin{array}{@{\,}c@{\,}}
u_{2}w^{8} \\
u_{3}w^{10} \\
\zeta w^{6} \\
 w^{7}
\end{array}
\right).
\label{4-29b}
\end{equation}
The coordinate transformation between the original coordinates and $(v_1,u_2,u_3,w, \zeta)$ is given by
\begin{equation}
\left\{ \begin{array}{l}
\displaystyle x_1 = w^{-2}  \\
\displaystyle x_2 = \left( \frac{2\sqrt{3}}{3} + v_1 \right) w^{-3},  \\
\displaystyle x_3 = \left( 2 + 3\sqrt{3}v_1 + u_{2}w^{8} - \frac{5\sqrt{3}}{8}w^7 + \frac{3}{8} zw^6 \right)w^{-4}, \\
\displaystyle x_4 =\left( \frac{8\sqrt{3}}{3} + 25v_1  + \frac{5\sqrt{3}}{3}u_{2}w^{8} + u_{3}w^{10}
                      + \frac{7\sqrt{3}}{8} zw^6 - \frac{27}{8} w^7 \right)w^{-5}, \\
\displaystyle z = \zeta.
\end{array} \right.
\label{4-30b}
\end{equation}
The equations of $v_1, u_2, u_3, w$ are
\begin{equation}
\left\{ \begin{array}{l}
\displaystyle \frac{dv_1}{dz} = \frac{\sqrt{3}v_1}{w}- \frac{3v_1^2}{2w}- \frac{5\sqrt{3}}{8}w^6 + u_2w^7 + \frac{3}{8}zw^5, \\[0.2cm]
\displaystyle \frac{du_2}{dz} =-\frac{3\sqrt{3}v_1^2}{2w^9} - \frac{15\sqrt{3}v_1}{16 w^2} + \frac{2v_1u_2}{w} + u_3w + \frac{3u_1z}{8 w^3}, \\[0.2cm]
\displaystyle \frac{du_3}{dz} =-\frac{15 v_1^2}{2w^{11}} + \frac{21 v_1}{16w^4} - \frac{5\sqrt{3}v_1u_2}{6} + \frac{5 v_1u_3}{2w} - \frac{3\sqrt{3}v_1z}{16 w^5}, \\[0.2cm]
\displaystyle \frac{dw}{dz} = -\frac{\sqrt{3}}{3} - \frac{1}{2}v_1. \\[0.2cm]
\end{array} \right.
\end{equation}
Although the right hand sides are not holomorphic at $w=0$,
they are holomorphic on the unstable manifold $v_1 = \varphi (u_2w^8, u_3w^{10}, zw^6, w^7) \sim O(w^6)$.
Any integral curves of the vector field outside the unstable manifold approach to the other fixed point (I).

\textbf{Step 5.} The $\Z_2$ action (\ref{Z2}) induces the $\Z_2$ action on the $(v_1, u_2, u_3, z, w)$ coordinates given by
\begin{eqnarray}
 (v_1, u_2, u_3, z, w) &\mapsto & ( -v_1 - \frac{4\sqrt{3}}{3}, u_2 - \frac{5\sqrt{3}}{4w} + \frac{6\sqrt{3}v_1}{w^8} + \frac{12}{w^8}, \\
& & -u_3 - \frac{10\sqrt{3}u_2}{3w^2} + \frac{25}{4w^3} - \frac{7\sqrt{3}z}{4 w^4}+ \frac{8\sqrt{3}}{w^{10}} - \frac{30 v_1}{w^{10}},z , -w). \qquad
\end{eqnarray}
If we divide $\C_2^4 = \{ (v_1, u_2, u_3, w)\}$ by this action, (\ref{4-30b}) becomes a one-to-one rational transformation
which defines a smooth algebraic variety $\C_0^4 \cup \C^4_2 /\Z_2$.
Therefore, $\mathcal{M}(z) = \C_0^4 \cup \C^4_1 /\Z_2 \cup \C^4_2 /\Z_2$ defined by the transformations (\ref{4-30}) and (\ref{4-30b})
gives the space of initial conditions for the system (\ref{4-14}).
\\


\textbf{Acknowledgements.}

The author would like to thank Professor Yasuhiko Yamada for useful comments.
This work was supported by Grant-in-Aid for Young Scientists (B), No.25800081 from MEXT Japan.


\appendix

\section{Normal form theory}

In this Appendix, we give a brief review of the normal form theory.
See \cite{Cho} for more detail.

Let us consider a holomorphic vector field
\begin{equation}
\frac{dx}{dt} = Ax + f(x), \quad x\in \C^m,
\label{A1}
\end{equation}
defined near the origin, where $A$ is an $m\times m$ matrix and $f\sim O(|| x ||^2)$ denotes the nonlinearity.
We assume that $A = \mathrm{diag} (\lambda _1,\cdots ,\lambda _m)$ is a diagonal matrix.
If there exist $j$ and non-negative integers $(n_1, \cdots ,n_m)$ such that $n_1 + \cdots +n_m \geq 2$ and
\begin{equation}
\lambda _1n_1 + \cdots + \lambda _mn_m = \lambda _j,
\end{equation}
then, the monomial vector field $x_1^{n_1}\cdots x_m^{n_m} \bm{e}_j$ is called the resonance term.
A normal form of (\ref{A1}) up to the order $N$ is given as follows.
\\[0.2cm]
\textbf{Proposition \thedef.} For any integer $N \geq 2$, there exists a polynomial transformation
$x \mapsto y$ of degree $N$ such that (\ref{A1}) is transformed into the system
\begin{equation}
\frac{dy}{dt} = Ay + g_1(y) + g_2(y),
\end{equation}
where $g_1$ consists only of resonance terms up to degree $N$, and $g_2 \sim O(|| x ||^{N+1})$.
\\[-0.2cm]

We need the following assumption for the convergence as $N\to \infty$.
\\[0.2cm]
\textbf{(P)} The convex hull of eigenvalues $\{ \lambda _1, \cdots ,\lambda _m\}$ in $\C$ does not include the origin.
\\[0.2cm]
In this case, the number of resonance terms is finite.
\\[0.2cm]
\textbf{Proposition \thedef.} Under the assumption (P), there exists a local analytic transformation
$x \mapsto y$ such that (\ref{A1}) is transformed into the system
\begin{equation}
\frac{dy}{dt} = Ay + g_1(y),
\end{equation}
where $g_1$ consists only of resonance terms.
\\[-0.2cm]

(\ref{A1}) has a \textit{formal} series solution of the form
\begin{equation}
x(t) = P(\alpha _1 e^{\lambda _1t}, \cdots ,\alpha _m e^{\lambda _mt}),
\end{equation}
where $P$ is a formal power series in the arguments, whose coefficients are polynomials in $t$.
$\alpha _1, \cdots , \alpha _m$ are arbitrary constants. 
The next proposition is well known in perturbation theory.
\\[0.2cm]
\textbf{Proposition \thedef.} $P$ is a \textit{convergent} power series, whose coefficients are independent of $t$, 
if and only if
\\[0.2cm]
(i) $A$ is semi-simple, and \\
(ii) (\ref{A1}) is linearized by a local analytic transformation.

In particular, the condition (ii) is satisfied if (P) is satisfied and $f(x)$ does not include resonance terms.
There are examples that (P) is not satisfied while (\ref{A1}) can be linearized (Siegel's theorem).
In the proof of Prop.3.5, the system (\ref{3-8}) satisfies (P) because the eigenvalues have positive real parts.


\section{Proofs of Proposition 2.7 and 2.8}

Consider the series solution
\begin{equation}
x_i(z) = c_i(z-z_0)^{-q_i} + a_{i,1} (z-z_0)^{-q_i+1} + \cdots 
 = c_iT^{-q_i}(1 + o(T)),
\label{B1}
\end{equation}
where $T = z-z_0$.
Without loss of generality, we suppose that 
\begin{equation}
\frac{q_1}{p_1} = \cdots  = \frac{q_M}{p_M} > \frac{q_{M+1}}{p_{M+1}} \geq \cdots \geq \frac{q_m}{p_m}.
\label{B2}
\end{equation}
Substituting (\ref{B1}) into $f_i^A$, we have
\begin{eqnarray*}
& & f_i^A (x_1(z) , \cdots ,x_m(z)) \\
&=& f_i^A (c_1 T^{-q_1}, \cdots ,c_mT^{-q_m}) \cdot (1 + o(T)) \\
&=& f_i^A (c_1 T^{\frac{q_M}{p_M}p_1 - q_1} T^{-\frac{q_M}{p_M}p_1}, \cdots 
            ,c_m T^{\frac{q_M}{p_M}p_m - q_m} T^{-\frac{q_M}{p_M}p_m} )  \cdot (1 + o(T)) \\
&=& T^{-\frac{q_M}{p_M}(1 + p_i)} f_i^A (c_1 T^{\frac{q_M}{p_M}p_1 - q_1},  \cdots 
                                         , c_m T^{\frac{q_M}{p_M}p_m - q_m} )  \cdot (1 + o(T)).
\end{eqnarray*}
(\ref{B2}) gives $q_Mp_i/p_M - q_i \geq 0$ for any $i$ and $q_Mp_i/p_M - q_i = 0$ for $i=1,\cdots ,M$.
Thus we obtain
\begin{eqnarray*}
f_i^A (x_1(z) , \cdots ,x_m(z)) = T^{-\frac{q_M}{p_M}(1 + p_i)} f_i^A (c_1 , \cdots ,c_M, 0, \cdots ,0) \cdot (1 + o(T)).
\end{eqnarray*}
Similarly, we can verify
\begin{eqnarray*}
& & f_i^N (x_1(z) , \cdots ,x_m(z),z) \sim o(T^{-\frac{q_M}{p_M}(1 + p_i)}), \\
& & g_i (x_1(z) , \cdots ,x_m(z),z) \sim o(T^{-\frac{q_M}{p_M}(1 + p_i)}),
\end{eqnarray*}
as $T\to 0$.
Hence, the system (\ref{2-1}) with (\ref{B1}) yields
\begin{equation}
-c_iq_i T^{-q_i-1} (1 + o(T)) = T^{-\frac{q_M}{p_M}(1 + p_i)} f_i^A (c_1 , \cdots ,c_M, 0, \cdots ,0) \cdot (1 + o(T)).
\label{B3}
\end{equation}
We will compare the orders of a pole in both sides; $q_i+1$ and $q_M(1+p_i)/p_M$.
\\[0.2cm]
\textbf{Proof of Prop.2.8.}
Assume that $q_i > p_i$ for some $i$.
We can assume without loss of generality that $q_i > p_i$ for any $i=1, \cdots ,m$.
Then, (\ref{B2}) shows
\begin{equation}
(q_i+1) - \frac{q_M}{p_M}(1+p_i) \leq (q_i+1) - \frac{q_i}{p_i}(1+p_i) = 1-\frac{q_i}{p_i} < 0.
\label{B4}
\end{equation}
This proves
\begin{eqnarray*}
f_i^A(c_1, \cdots , c_M, 0 ,\cdots ,0) = 0
\end{eqnarray*}
for $i=1, \cdots ,m$.
Then, the condition (S) gives $c_1 = \cdots = c_M = 0$.
We repeat this procedure by replacing $q_i-1$ by $q_i$ ($q_i-1 := q_i$) for $i=1,\cdots ,M$ and 
rearranging the order of $q_1, q_2, \cdots ,q_m$ so that (\ref{B2}) holds for some $M$.
If $q_i > p_i$, the inequality (\ref{B4}) again holds.
If $q_i = p_i$ for some $i$, then
\begin{eqnarray*}
(q_i+1) - \frac{q_M}{p_M}(1+p_i) < (q_i+1) - \frac{q_i}{p_i}(1+p_i) = 0.
\end{eqnarray*}
Since the inequality (\ref{B4}) holds for any $i$, we have $c_1 = \cdots =c_M = 0$.

By repeating this procedure, at least $q_1$ decreases by $1$ at each step.
This algorithm stops when $q_i = p_i$ for any $i$ and it completes the proof. $\Box$
\\[0.2cm]
\textbf{Proof of Prop.2.7.}
Suppose $0\leq q_i < p_i$ and $c_i\neq 0$ for $i=1, \cdots ,m$.
For $i = 1, \cdots ,M$, we have
\begin{eqnarray*}
(q_i+1) - \frac{q_M}{p_M}(1+p_i) = (q_i +1) - \frac{q_i}{p_i}(1 + p_i) = 1-\frac{q_i}{p_i} > 0.
\end{eqnarray*}
Thus, Eq.(\ref{B3}) gives $c_iq_i = 0$ for $i=1,\cdots ,M$.
By the assumption, we obtain $q_i=0$ for $i=1,\cdots ,M$.
By repeating this procedure as before, we can prove that $q_i = 0$ for any $i$. $\Box$


\end{document}